\documentclass[a4paper,11pt]{article}
\usepackage{textcomp}
\usepackage[left=20mm, top=20mm, right=20mm, bottom=25mm]{geometry}
\usepackage{graphicx}
\usepackage{url}
\usepackage{enumitem}

\usepackage{breakcites}
 
\usepackage[usenames,dvipsnames]{color}
\usepackage{amssymb}
\usepackage{amsmath}
\usepackage{array}
\usepackage{tabularx}
\usepackage{hyperref,wasysym}
\usepackage{scrextend}
\usepackage{multirow}
\usepackage{smartdiagram}
\usepackage{metalogo}
\usepackage{tikz}
\usepackage{caption}
\usetikzlibrary{intersections, backgrounds}
\usetikzlibrary{calc, positioning}
\usetikzlibrary{shapes.arrows, fadings, shadows}
\usetikzlibrary{shapes.geometric}
\usepackage[utf8]{inputenc}
\usepackage[T1]{fontenc}
\usepackage{booktabs}
\usepackage{multicol}
\usepackage{adjustbox}
\usepackage[leftmargin=0em]{quoting}
\usepackage[round]{natbib} 
\definecolor{azure}{rgb}{0.2, 0.4, 1.0}
\usepackage{color}
\usepackage{xargs}
\usepackage[colorinlistoftodos,prependcaption,textsize=normalsize]{todonotes}
\newcommandx{\unsure}[2][1=]{\todo[linecolor=red,backgroundcolor=red!25,bordercolor=red,#1]{#2}}
\newcommandx{\change}[2][1=]{\todo[linecolor=blue,backgroundcolor=blue!25,bordercolor=blue!50,#1]{#2}}
\newcommandx{\info}[2][1=]{\todo[linecolor=OliveGreen,backgroundcolor=OliveGreen!25,bordercolor=OliveGreen,#1]{#2}}
\newcommandx{\improvement}[2][1=]{\todo[linecolor=Plum,backgroundcolor=Plum!25,bordercolor=Plum,#1]{#2}}
\newcommandx{\thiswillnotshow}[2][1=]{\todo[disable,#1]{#2}}

\tikzfading[name=arrowfading, top color=transparent!0, bottom color=transparent!95]
\tikzset{arrowfill/.style={top color=OrangeRed!20, bottom color=Red, general shadow={fill=black, shadow yshift=-0.8ex, path fading=arrowfading}}}
\tikzset{arrowstyle/.style={draw=FireBrick,arrowfill, single arrow,minimum height=#1, single arrow,
single arrow head extend=.4cm,}}

\title{\vspace*{-1.0cm}\textbf{English-Medium Instruction Calculus: Is flipping helpful?}}
\author{\color{black}\normalsize N. Karjanto{\footnotesize $^{1}$}\thanks{Corresponding author} \; and L. Simon{\footnotesize $^{2}$}\\ 
\color{black}\small Department of Mathematics, University College\\
\color{black}\small Sungkyunkwan University, Natural Science Campus\\
\color{black}\small 2066 Seobu-ro, Suwon-si, Jangan-gu, 16419\\
\color{black}\small Gyeonggi-do, Republic of Korea\\
\color{black}\small {\footnotesize $^{1}$}\url{natanael@skku.edu}, {\footnotesize $^{2}$}\url{loisssimon@skku.edu}
}
\date{\scriptsize Updated \today}

\begin{document}
\maketitle
\begin{abstract}
This paper addresses some of the experiences encountered and some of the strategies utilized while using flipped classroom pedagogy in an undergraduate Single Variable Calculus (SVC) course in a Confucian Heritage Culture (CHC) environment. Due to issues related to instructor-student communication and students' learning styles, a theoretical framework is built upon three components for learning environments: inverted Bloom's taxonomy for educational learning objectives, English-medium instruction (EMI) and integrating technology in the classroom. In this study, four different types of instruction were designed; three classes were flipped while the fourth acts as a control. Type A means entirely flipped with instructor-created videos and integrating CAS {\sl Maxima}. Type B means entirely flipped with third-party videos. Type C is similar to Type A, but instead of the entire course being flipped, a single topic is flipped. Type D uses the traditional lecture style.  By considering exam scores and students' feedback, both quantitative and qualitative aspects of the flipped pedagogy were investigated. Quantitative analysis indicates a statistically significant difference in the exam scores between instruction Type A and Type C as determined by one-way ANOVA ($F(3,306) = 2.67$, $p = 0.0477$) with small practical significance of the effect size ($\eta^2 = 0.0255$). Qualitative findings suggest that although students' engagement and communication with instructors have been enhanced, some challenges related to language and culture, along with minimizing competition and learning about a new technological tool remain to persist.

\end{abstract}

\section{Introduction}

The term `flipped classroom' has been popularized in the past decade among many educators who have sought to introduce innovative teaching methods in their classes. The term was coined and popularized by Jonathan Bergmann and Aaron Sams in 2007~\citep{bergmann2012flip}. A few examples of similar terminology used in the literature are `inverted classroom'~\citep{lage2000inverting, talbert2014inverting}, `learn before lecture'~\citep{moravec2010learn}, `flipped class'~\citep{fulton2012upside}, `flipped learning'~\citep{marshal2014,dodds2015evidence}, and `flipping pedagogy'~\citep{mcgivney2013flipping}; there are many more. In the context of this article, they refer to a unified idea, i.e. `lecture at home' and `homework at school'. For the remainder of this paper, we use `flipped classroom'.

The flipped classroom method is not new, but the approach has increased in popularity due to advancements in technology~\citep{fitzpatrick2012}. First, students gain exposure to new concepts outside the class, through reading assigned materials, watching recorded video lectures, and/or listening to podcasts~\citep{berrett2012flipping}. Second, class time is dedicated to addressing the more difficult task of assimilating that knowledge by engaging in activities such as problem-solving, critical thinking, experiments, discussion or debates~\citep{brame2013flipping}. The governing board and key leaders of the Flipped Learning Network describe this as a pedagogical approach where direct instruction is relocated from a collective learning environment to an individual learning environment. The former is then transformed into a dynamic, interactive learning environment where educators guide students in applying concepts and engaging creatively with the materials~\citep{fln2016}.

Although the flipped classroom has increased in popularity among educators it is not always as popular with learners.  Some reasons for unpopularity among students include requiring students to actively collect information prior to attending class, along with requiring active student participation during class time. Other concerns include feelings of abandonment that result when students are encouraged to practice self-study in contrast to the students' cultural learning environment where lecturing and rote learning are normative~\citep{talbert2012inverted}. 

Nevertheless, the flipped classroom has great potential to foster interactive and effective learning activities for all learners by incorporating peer instruction and discussion. This potential includes improving the traditional lecture~\citep{berrett2012flipping}, boosting homework completion rates~\citep{moore2014fostering}, accommodating various learning styles~\citep{bergmann2012flip,bergmann2014flipped}, altering the focus from information transfer to information assimilation~\citep{mazur2009farewell}, as well as the positive effect of improving students' academic performance~\citep{crouch2001peer,smith2009peer,deslauriers2011improved}. Additionally, the latter provides supporting evidence in science-related subjects, in connection to peer instruction and discussion.

In the context of active learning, the largest and most comprehensive metaanalysis study on active learning, student performance and failure rate in the STEM field has been reviewed by~\cite{freeman2014active}. Their findings indicate that on average, student academic performance increased by a standard deviation with an effect size of 0.47 and examination scores improved by 6\% in active learning sessions. Alternatively, students in the traditional classroom were more likely to fail than their peers in an active learning classroom, with the risk ratio of 1.5 and an effect size of 1.95. In light of these findings, we consider that a classroom with active learning components and interactive engagement activities may prove to be more effective than the traditional classroom. Indeed, by adopting flipped classroom pedagogy within our SVC courses, we were hoping to promote active learning along with enhanced communication in the classroom.
 
In this study, we are interested in investigating whether non-native English speakers, particularly those in the CHC context, perform better in the flipped classroom or the traditional classroom.
Hence, we attempt to answer the following research questions:
\begin{itemize}[leftmargin=1em]
\item (Quantitative) Will students in the flipped classroom academically outperform their traditional classroom peers? 
\item (Quantitative) Is there a statistically significant difference in academic performance among students in a flipped class with CAS vs. those in a flipped class without CAS?
\item (Quantitative) Is there a statistically significant difference in terms of academic performance among students viewing instructor-created videos vs. students viewing third-party created videos?
\item (Qualitative) Has our flipped classroom pedagogy successfully bridged the educator-learner cultural gap by improving interaction and communication between students and instructors in an EMI--CHC context?
\end{itemize}

This article comprises seven sections. Following the Introduction, Section~\ref{litera} provides a literature review on flipped classroom, in mathematics and various fields. Section~\ref{TF} outlines our study's theoretical framework. The original, revised and the inverted Bloom's taxonomies are discussed. It also covers issues addressing English as the medium of instruction. Technology usage, particularly in mathematics education, is also highlighted. Section~\ref{CA} describes flipped classroom pedagogy in SVC course at Sungkyunkwan University (SKKU). This section elucidates course logistics; course topics, the adopted textbook, recorded videos, technology usage in mathematics education through CAS {\sl Maxima} and an intranet learning management system (LMS) \textsl{icampus}. As we proceed, Section~\ref{method} comprehensively deals with the research methodology employed. This includes, descriptive participant statistics, measurements, data collection and statistical analysis of the quantitative aspect of the study. In Section~\ref{finding}, we report quantitative and qualitative findings resulting from the implementation of flipped classroom pedagogy for SVC in an EMI--CHC context. Finally, Section~\ref{discussion} provides some of the study's drawbacks and limitations, concluding remarks, and future implications of our findings.

\section{Literature Review} \label{litera}

Flipped classroom pedagogy has been implemented in various subjects at the secondary and tertiary levels. In undergraduate mathematics, this includes College Algebra, Calculus and Linear Algebra. Other fields implementing flipped classrooms include Economics, Engineering, Computer Science, Medicine and Biological Science. The following provides a non-exhaustive literature overview of the flipped classroom encompassing various levels. In K-12 education, a flipped classroom success story in PreCalculus, Calculus and (Accelerated) Algebra and Geometry has been documented by~\cite{fulton2012upside}. At the tertiary level, a quasi-experimental quantitative research study in College Algebra comparing the flipped classroom and traditional lecture and its effect on student achievement as measured through common assessments has been independently investigated by~\cite{overmyer2014flipped}, \cite{van2015adventures} and~\cite{acelajado2017flipped}. 

There are numerous studies in Calculus at the tertiary level. \cite{mcgivney2013flipping} have discussed students' perceptions of the flipped classroom on definite integral applications while \cite{sahin2015flipping} have compared student performance between flipped and traditional classes. Analogously, \cite{ziegelmeier2015flipped}~have explored Multivariable Calculus. On the other hand, \cite{jungic15} investigated flipped classroom with large class sizes and acknowledged that although a significant amount of time and resources is invested for the preparation, it is beneficial and gratifying. \cite{maciejewski2016flipping}~discussed students' mathematical understanding and achievement; flipped classroom students on average outperformed their peers from the traditional classroom by 8\% in academic performance. Further analysis indicates that students with high mathematical ability but limited exposure to Calculus benefit most from the experiment. Interestingly, regarding students' attitude toward mathematics, the Calculus flipped classroom works best with students who have previous experience and a solid Calculus background, as argued by~\cite{sonnert2015impact}. 

In Linear Algebra, \cite{talbert2014inverting} proposed a number of flipped classroom designs, i.e. flipping a single-topic, a series of workshop, and flipping an entire course. 
\cite{murphy2016student} conducted a further study in flipped classroom Linear Algebra and confirmed that students' attitude is positive, as indicated by positive enjoyment, more confidence to study independently and to exhibit better retention rate of the materials. Recently, \cite{novak2017flip} reported students' perspective of a flipped classroom in Linear Algebra on the topic of calculating matrix determinants. Greater engagement and increased understanding of the material covered are reported.

There are various articles including and excluding STEM subjects; for further insight, we provide a small sample, Actuarial Science~\citep{butt2014student}, Business~\citep{findlay2014evaluation}, Biology~\citep{moravec2010learn}, Chemistry~\citep{fautch2015flipped}, Computer Science~\citep{davies2013flipping,mok2014teaching,dodds2015evidence}, Economics~\citep{lage2000inverting}, Engineering and Mathematics~\citep{lape2014probing}, Physics~\citep{kettle2013flipped}, Pharmacy~\citep{pierce2012vodcasts}, Science Education~\citep{herreid2013case}, and Statistics~\citep{wilson2013flipped}. 

Another comparative study of the flipped classroom in three different fields of Engineering, Social Studies and Humanities has been explored by~\cite{kim2014experience}. In addition to examining multi-disciplinary applications of the flipped classroom, they also proposed a design framework and nine design principles for the pedagogical approach. A recent study from the same field on the effectiveness of the flipped classroom in comparison to blended learning, traditional learning and e-learning has been explored by~\cite{thai2017impact}. A comprehensive survey of flipped classroom research has been provided by~\cite{bishop2013flipped}. They reported that student perceptions are generally positive. \cite{chen2016exploring} explored the relationship between student perceptions, gender difference and learning outcomes. 

Another literature survey on the flipped classroom, particularly in STEM subjects, has been explained by~\cite{dodds2015evidence}. Additionally, the effectiveness of a flipped classroom model on student engagement and achievement, as well as the affordance of this model versus a traditional one, has been explored by~\cite{bormann2014affordances}. 
\cite{loucky2017flipped} recently published and explored the flipped classroom with educational and digital technologies for effective second language learning curricula. Moreover, \cite{towey2015lessons} and \cite{ryan2013socrates} focused on the challenges related to implementing the flipped classroom in the EMI--CHC context, for example the language barrier and cultural differences, particularly the passive-receptive learning style.

\section{Theoretical Framework} \label{TF}

The theoretical framework for this study is built from an intersection of three components, inverted Bloom's taxonomy of educational learning objectives, EMI for non-native English speakers in the CHC context and the adaptation of technology for mathematical instruction. The latter uses instructor-created and {\sl Khan Academy}-created video recordings, intranet LMS {\sl icampus}, and a Computer Algebra System (CAS) {\sl Maxima}. Figure~\ref{venn2} summarizes the theoretical framework for our flipped classroom pedagogy. In what follows, we will discuss each component.
\def\firstcircle{(3cm,1.7cm) circle (1.3cm)}
\def\secondcircle{(3cm,-1.7cm) circle (1.3cm)}
\def\thirdcircle{(0cm,0cm) circle (1.3cm)}
\begin{figure}[h]
\begin{center}
\begin{tikzpicture}
\begin{scope}[fill opacity = 0.7]
  \fill[red!90!orange] \firstcircle;
  \fill[green!70!black] \secondcircle;
  \fill[blue!50!cyan] \thirdcircle;
\end{scope}
  \node at (3.2cm,1.7cm)  [opacity = 0.8, text width=1.7cm]{\small \sffamily Inverted Bloom's taxonomy};
  \node at (2.9cm,-1.7cm) [opacity = 0.8, text width=1.7cm]{\small \sffamily EMI~in~CHC};
  \node at (0.2cm,0cm)    [opacity = 0.8, text width=2.0cm]{\small \sffamily Technology adaptation {\hspace*{0.2cm} (CAS)}};
\node at (7cm,0cm) [
  right color=blue!70!red,
  left color=red!60!blue,
  single arrow,
  minimum height=2.5cm,
  minimum width=2cm,
  shading angle=90,
  rotate=0
]{}; 
\end{tikzpicture}
\smartdiagramset{set color list={red!90!orange, green!70!black, blue!50!cyan}, 
   bubble center node size = 3cm, 
   bubble center node font = \normalfont \sffamily, 
   bubble center node color = red!30!yellow, 
   distance center/other bubbles = 0.6cm,
   distance text center bubble = 0.5cm,
   bubble node font = \small \sffamily,
   bubble fill opacity = 0.7}
\smartdiagram[bubble diagram]{Flipped\\ classroom, Inverted\\ Bloom's\\ taxonomy, EMI in CHC, Technology\\ adaptation\\ (CAS)}
\caption{\small A theoretical framework for flipped classroom pedagogy implemented in this study is built from an intersection of three components: inverted Bloom's taxonomy of educational learning objectives, EMI in CHC, and technology adaptation, particularly using CAS.} \label{venn2}
\end{center}
\vspace*{-0.75cm}
\end{figure}
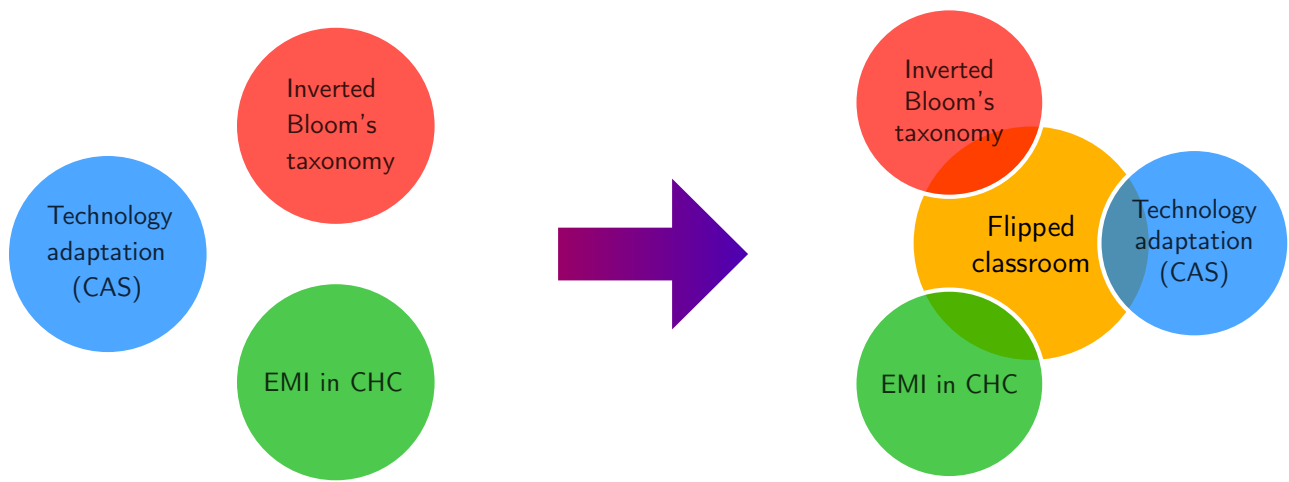

\subsection{Inverted Bloom's taxonomy}

\begin{figure}[h]
\begin{center} 
\begin{tikzpicture}[rounded corners, scale = 0.8, every node/.style={scale=0.8}]
\draw[fill = green!75!blue, draw = blue!70!green] (0,0)--(8.2,0)--(7.6,0.9)--(0.6,0.9)--cycle node at (4.0,0.4) {\sffamily Remember};
\draw[fill = green!60!blue, draw = blue!70!green] (0.7,1)--(7.5,1)--(6.9,1.9)--(1.3,1.9)--cycle node at (4.0,1.4) {\sffamily Understand};
\draw[fill = green!50!blue, draw = blue!70!green] (1.4,2)--(6.8,2)--(6.2,2.9)--(2.0,2.9)--cycle node at (4.0,2.4) {\sffamily Apply};
\draw[fill = blue!60!green, draw = blue!70!green] (2.1,3)--(6.1,3)--(5.5,3.9)--(2.7,3.9)--cycle node at (4.0,3.4) {\sffamily Analyze};
\draw[fill = blue!70!green, draw = blue!70!green] (2.8,4)--(5.4,4)--(4.8,4.9)--(3.4,4.9)--cycle node at (4.1,4.4) {\sffamily \color{white} Evaluate};
\draw[fill = blue!80!green, draw = blue!70!green] (3.5,5)--(4.7,5)--(4.1,5.9)--cycle node at (5.0,5.4) {\sffamily Create};
\node at (9,2.7) [draw = blue!70!green,
    right color= blue!80!green,  
    left color = green!80!blue,
    single arrow,
    minimum height=5.9cm,
    minimum width=1.8cm,
    shading angle=180,
    rotate=90
] {\sffamily Cognitive Process Dimension};
\end{tikzpicture}
\includegraphics[width = 0.45\textwidth]{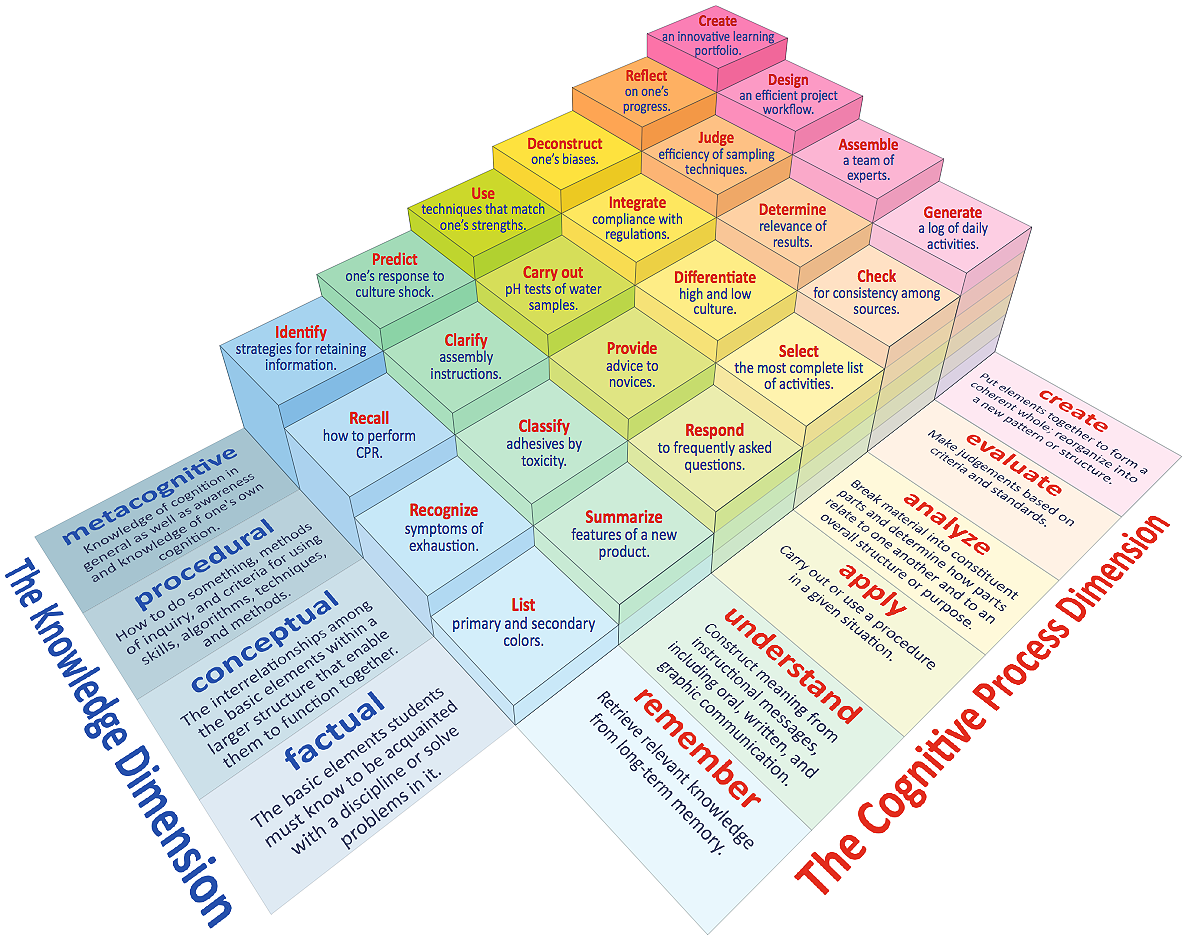}
\end{center}
\caption{A pyramid representation of the (revised) Bloom's taxonomy of educational learning objectives covering cognitive process dimension (left panel) and a three-dimensional representation of the taxonomy which includes both cognitive process and knowledge dimensions (right panel). The level of difficulty of acquiring learning objectives increases from bottom to top.}
\label{bloom1}
\end{figure}
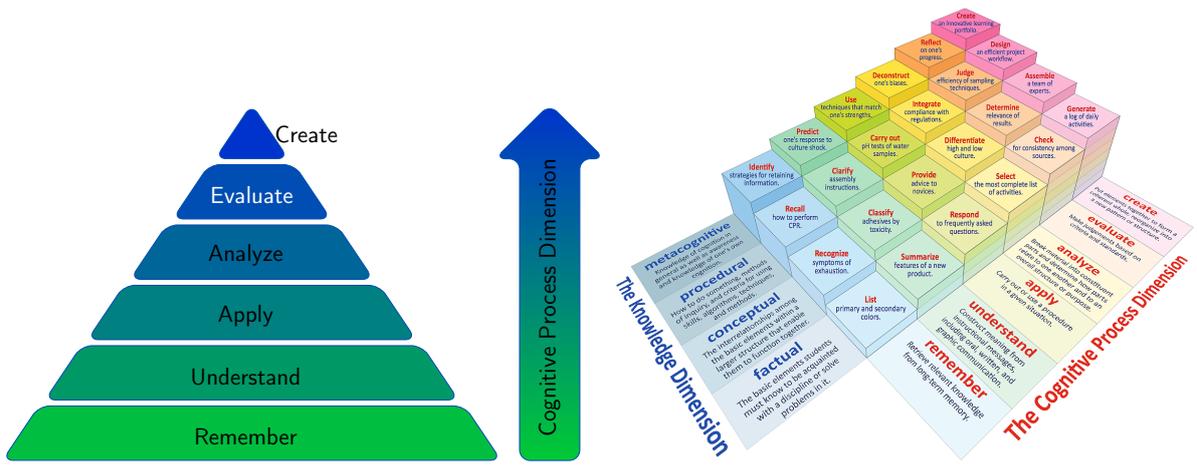

The first essential component we will discuss is in the context of Bloom's taxonomy of educational learning objectives~\citep{bloom1964taxonomy}; the revised version has been discussed by~\cite{anderson2001taxonomy} and~\cite{krathwohl2002revision}. In Figure~\ref{bloom1}, a pyramid representation of the revised version covering the cognitive process dimension and a three-dimensional representation of the taxonomy covering both the cognitive process and knowledge dimensions are displayed. In the cognitive process dimension, the lower-order thinking skills of `Remember' and `Understand' are easier to acquire than the higher-order thinking skills of `Evaluate' and `Create'. Similarly, in the knowledge dimension, the `Factual' and `Conceptual' aspects are easier to master than the `Procedural' and `Metacognitive' aspects.

In the traditional classroom, easier processes are generally acquired in the classroom while more difficult processes are delegated to the students outside of the classroom, usually via assignments. Hence, when students need assistance from their instructors the most, while acquiring higher-order thinking skills, they are largely abandoned. For many disciplines with strong relationships among topics, including Mathematics and other science-related fields, the task of acquiring higher-order learning objectives may be a source of frustration to those students who need hands-on assistance and guidance.

In the flipped classroom, we reversed the time and place for these learning activities. First, we delegate the easier learning objectives and acquisition of the lower-order thinking skills to students outside of the classroom. Second, we dedicate the in-class time to learning activities that aid in acquiring the higher-order thinking skills. Therefore, the skills `Remember' and `Understand' can be acquired before class, the skills `Apply' and `Analyze' can be developed during the class and the skills `Evaluate' and `Create' can be acquired after class once the students gain confidence and are equipped during in-class sessions. For in-class activities, we strive to create a learning environment where students could receive assistance and overcome difficulties in acquiring more difficult cognitive skills. This assistance can be provided through the instructor's guidance as well as through peer collaborative learning~\citep{king1993sage,saulnier2008sage}. By implementing this pedagogy, our aim is not only to assist frustrated students who struggle with homework but also to cultivate self-discipline in the entire class through preparation before the class and interactive engagement during the class. 
\begin{figure}[htp]
\begin{center}
\begin{tikzpicture}[rounded corners, scale = 0.8, every node/.style={scale=0.8}]
\draw[fill = green!75!blue, draw = blue!70!green] (0,6)  --(8.2,6)--(7.6,5.1)--(0.6,5.1)--cycle node at (4.1,5.6) {\sffamily Remember};
\draw[fill = green!60!blue, draw = blue!70!green] (0.7,5)--(7.5,5)--(6.9,4.1)--(1.3,4.1)--cycle node at (4.1,4.6) {\sffamily Understand};
\draw[fill = green!50!blue, draw = blue!70!green] (1.4,4)--(6.8,4)--(6.2,3.1)--(2.0,3.1)--cycle node at (4.0,3.5) {\sffamily Apply};
\draw[fill = blue!60!green, draw = blue!70!green] (2.1,3)--(6.1,3)--(5.5,2.1)--(2.7,2.1)--cycle node at (4.0,2.5) {\sffamily Analyze};
\draw[fill = blue!70!green, draw = blue!70!green] (2.8,2)--(5.4,2)--(4.8,1.1)--(3.4,1.1)--cycle node at (4.1,1.5) {\sffamily \color{white} Evaluate};
\draw[fill = blue!80!green, draw = blue!70!green] (3.5,1)--(4.7,1)--(4.1,0.1)--cycle node at (5.0,0.5) {\sffamily Create};
\node at (9,3.3) [draw = blue!70!green,
    right color= green!80!blue,  
    left color = blue!80!green,
    single arrow,
    minimum height=5.9cm,
    minimum width=1.8cm,
    shading angle=180,
    rotate=-90
] {\sffamily Cognitive Process Dimension};
\end{tikzpicture}
\begin{tikzpicture}[rounded corners, scale = 0.8, every node/.style={scale=0.75}]
\draw[fill = blue!80!green, draw = blue!70!green] (0,0) rectangle node{\sffamily \color{white} Create} ++(1.9,0.7);
\draw[fill = none, draw = blue!80!green, line width=0.7mm] (2,0) rectangle node{\sffamily Generate} ++(1.9,0.7);
\draw[fill = none, draw = blue!80!green, line width=0.8mm] (4.02,0) rectangle node{\sffamily Assemble} ++(1.88,0.7);
\draw[fill = none, draw = blue!80!green, line width=0.9mm] (6.05,0) rectangle node{\sffamily Design} ++(1.9,0.7);
\draw[fill = none, draw = blue!80!green, line width=1.0mm] (8.1,0) rectangle node{\sffamily Create} ++(2.45,0.7);
\draw[fill =  blue!70!green, draw =  blue!70!green] (0,1) rectangle node{\sffamily \color{white} Evaluate} ++(1.9,0.7);
\draw[fill =  none, draw =  blue!70!green, line width=0.6mm] (2,1) rectangle node{\sffamily Check} ++(1.9,0.7);
\draw[fill =  none, draw =  blue!70!green, line width=0.7mm] (4.01,1) rectangle node{\sffamily Determine} ++(1.89,0.7);
\draw[fill =  none, draw =  blue!70!green, line width=0.8mm] (6.05,1) rectangle node{\sffamily Judge} ++(1.9,0.7);
\draw[fill =  none, draw =  blue!70!green, line width=0.9mm] (8.1,1) rectangle node{\sffamily Reflect} ++(2.45,0.7);
\draw[fill = blue!60!green, draw =  blue!70!green] (0,2) rectangle node{\sffamily Analyze} ++(1.9,0.7);
\draw[fill = none, draw =  blue!60!green, line width=0.5mm] (2,2) rectangle node{\sffamily Select} ++(1.9,0.7);
\draw[fill = none, draw =  blue!60!green, line width=0.6mm] (4,2) rectangle node{\sffamily Differentiate} ++(1.95,0.7);
\draw[fill = none, draw =  blue!60!green, line width=0.7mm] (6.1,2) rectangle node{\sffamily Integrate} ++(1.85,0.7);
\draw[fill = none, draw =  blue!60!green, line width=0.8mm] (8.1,2) rectangle node{\sffamily Deconstruct} ++(2.45,0.7);
\draw[fill = green!50!blue, draw =  blue!70!green] (0,3) rectangle node{\sffamily Apply} ++(1.9,0.7);
\draw[fill = none, draw =  green!50!blue, line width=0.4mm] (2,3) rectangle node{\sffamily Respond} ++(1.9,0.7);
\draw[fill = none, draw =  green!50!blue, line width=0.5mm] (4,3) rectangle node{\sffamily Provide} ++(1.9,0.7);
\draw[fill = none, draw =  green!50!blue, line width=0.6mm] (6.05,3) rectangle node{\sffamily Carry out} ++(1.9,0.7);
\draw[fill = none, draw =  green!50!blue, line width=0.7mm] (8.1,3) rectangle node{\sffamily Use} ++(2.45,0.7);
\draw[fill =  green!60!blue, draw =  blue!70!green] (0,4) rectangle node{\sffamily Understand} ++(1.9,0.7);
\draw[fill =  none, draw =  green!60!blue, line width=0.3mm] (2,4) rectangle node{\sffamily Summarize} ++(1.9,0.7);
\draw[fill =  none, draw =  green!60!blue, line width=0.4mm] (4,4) rectangle node{\sffamily Classify} ++(1.9,0.7);
\draw[fill =  none, draw =  green!60!blue, line width=0.5mm] (6.05,4) rectangle node{\sffamily Clarify} ++(1.9,0.7);
\draw[fill =  none, draw =  green!60!blue, line width=0.6mm] (8.1,4) rectangle node{\sffamily Predict} ++(2.45,0.7);
\draw[fill =  green!75!blue, draw =  blue!70!green] (0,5) rectangle node{\sffamily Remember} ++(1.9,0.7);
\draw[fill =  none, draw =  green!75!blue, line width=0.2mm] (2,5) rectangle node{\sffamily List} ++(1.9,0.7);
\draw[fill =  none, draw =  green!75!blue, line width=0.3mm] (4,5) rectangle node{\sffamily Recognize} ++(1.9,0.7);
\draw[fill =  none, draw =  green!75!blue, line width=0.4mm] (6.05,5) rectangle node{\sffamily Recall} ++(1.9,0.7);
\draw[fill =  none, draw =  green!75!blue, line width=0.5mm] (8.1,5) rectangle node{\sffamily Identify} ++(2.45,0.7);
\draw[fill =  yellow!50, draw =  red!80!yellow] (2,6) rectangle node{\sffamily Factual} ++(1.9,0.7);
\draw[fill =  yellow!50!red, draw =  red!80!yellow] (4,6) rectangle node{\sffamily Conceptual} ++(1.9,0.7);
\draw[fill =  red!70!yellow, draw =  red!80!yellow] (6.05,6) rectangle node{\sffamily Procedural} ++(1.9,0.7);
\draw[fill =  red!90!yellow, draw =  red!80!yellow] (8.1,6) rectangle node{\sffamily Metacognitive} ++(2.5,0.7);
\node at (6,7.5) [
    right color=red!80!black,
    left color=orange!50!yellow,
    single arrow,
    minimum height=7.5cm,
    minimum width=1.5cm,
    shading angle=90,
] {\sffamily \large Knowledge Dimension};
\end{tikzpicture}
\end{center}
\caption{A reversed pyramid of Bloom's taxonomy for cognitive process dimension (left panel) and a two-dimensional representation of Bloom's taxonomy that includes both the cognitive process and knowledge dimensions. The level of difficulty in acquiring learning objectives rises from top to bottom and from left to right for the cognitive process and knowledge dimensions, respectively.} \label{bloom2}
\end{figure}
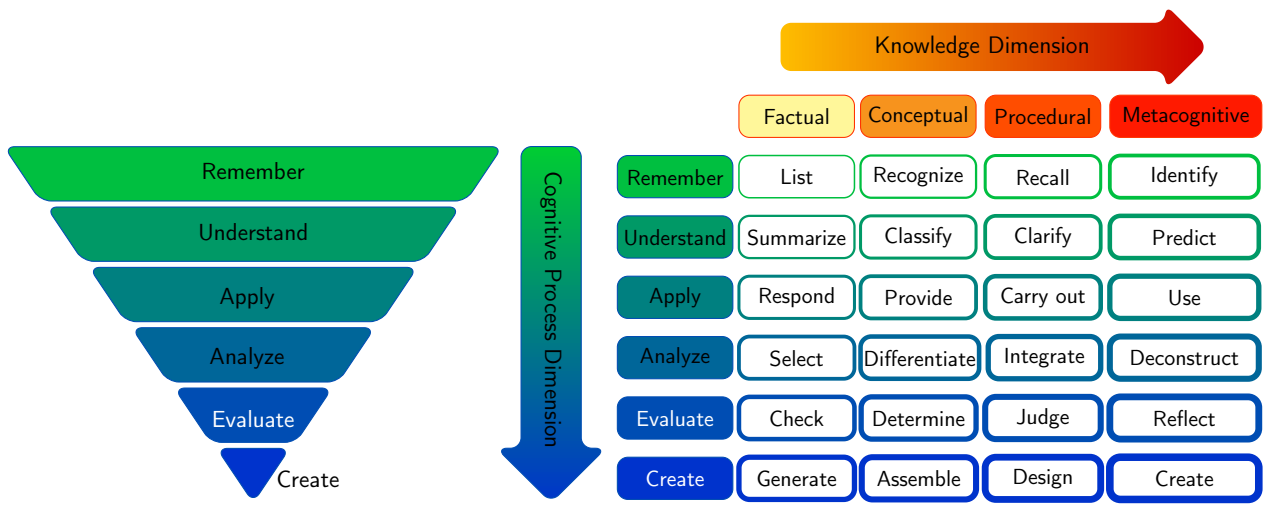

Although \cite{bergmann2014flipped} used the term `inverted Bloom's taxonomy', we remark that what has been actually inverted is not Bloom's taxonomy itself, but rather the time (before class vs. during class and during class vs. after class) and the location (in-class vs. out-of-class) for the learning activities. A learning process for flipped classroom pedagogy is described by a reversed pyramid of Bloom's taxonomy as shown in the left panel of Figure~\ref{bloom2}. The right panel of Figure~\ref{bloom2} shows a two-dimensional representation of Bloom's taxonomy which includes both the cognitive process and knowledge dimensions. The level of difficulty in acquiring learning objectives goes from top to bottom and from left to right for the cognitive process and knowledge dimensions, respectively. In addition, we also remark that in the presence of language and cultural barriers, acquiring `List', `Recognize', `Summarize', or `Classify' may be higher-order thinking skills for non-native English speakers accustomed to a quiet class culture and a passive-receptive learning environment. In the following subsection, we will address other EMI concerns in the CHC context.

\subsection{English-medium instruction (EMI)} \label{EMI}

At SKKU, SVC, titled Calculus~1 (course code GEDB001) is designated as an international language course. The chosen language of instruction is English and the majority of students are non-native English speakers. Generally speaking, it is difficult to learn another language. Many university-aged Korean students have had exposure to English from early childhood to late teen-hood. Their exposure is a consequence of the local school curriculum and also after-school programs and private academies, locally known as `hagwons'. Although, university students may have had long-term exposure to English, we observe that some students still have some difficulties while attending an EMI course. These difficulties are not as pronounced among international students who have had exposure to English. In this subsection, we discuss the second component of our theoretical framework, i.e. delivering a course while using English as the medium of instruction where English is not the first language for the majority of the students.

Many university programs around the world have adopted English as a medium of instruction as part of their programs.
This bold step has fostered and encouraged debates, challenges and controversies~\citep{doiz2012english,gundermann2014}.
EMI has been discussed for a number of European universities during the last decade, for example in Austria~\citep{tatzl2011english}, Denmark~\citep{werther2014using}, Germany~\citep{klippel2003new}, Italy~\citep{costa2013survey}, the Netherlands~\citep{vinke1998english,wilkinson2013english}, Spain~\citep{cots2013,doiz2011internationalisation}, Sweden~\citep{bjorkman2008so,bolton2012english} and Switzerland~\citep{murray2001dominance} as well.
Some challenges include, but are not limited to, improving language skills, instructor's competencies, difficulty in implementation and the fear of native language survival. Nevertheless, the literature also suggests that EMI delivers ample benefits: improving English listening and speaking skills~\citep{tatzl2011english}, increasing multilingual competency, allowing for the practice of a wider range of study skills and learning strategies~\citep{klippel2003new} and internationalization in higher education. Some institutions provide training for instructors to support EMI implementation~\citep{klaassen2001facing}. For a broader look at EMI in an European higher education context, readers are encouraged to consult~\cite{coleman2006english}.

In the Republic of Korea (South Korea), this practice is prevalent because	 it is related to the government's accreditation process for higher education institutions. This practice is also common because offering EMI courses may boost the university's national ranking and thus raises its popularity among potential students, local and international, as well as with grant sponsors and university donors~\citep{cho2012english}. \cite{byun2011english} examined the effectiveness of the EMI policy in the South Korean higher education system and discovered that the lack of language proficiency from both students and instructors is prevalent. \cite{kang2012english} provided evidence against the exclusive use of EMI in South Korean classrooms and highlighted the benefits of teaching using the mother tongue. \cite{kim2014emerging} stated that the students enrolled in EMI courses exhibit very low confidence in interactive activities. Language proficiency significantly affects academic performance. On the other hand, \cite{joe2013does} argued that for medical students, the medium of instruction had no effect on students' understanding of the lecture and students' proficiency was not correlated to their comprehension.

As the literature has indicated, limited language proficiency may hinder the teaching and learning process. We observed that there exists an additional factor to this hindrance, namely the students' culture. Although many parties express favorable attitudes toward EMI, some may not necessarily favor its practices, as recently revealed by~\cite{huang2018resistance}. The author found that most Taiwanese learners resisted curriculum, pedagogy, and context. Their resistance is strongly related to the CHC learning style. CHC students are often portrayed as quiet, passive, with reluctance in expressing opinions and are not accustomed to actively participating in the classroom, even though not all learners possess homogeneous learning styles and the criticism has been challenged by~\cite{kember2000misconceptions,kennedy2002learning,tran2013learning}. Nevertheless, from our teaching experiences with some of our students in the South Korean classrooms, some of the former stereotypes are often encountered.

In this study, by implementing flipped classroom pedagogy in our SVC courses, we hope to overcome some of the challenges students face when learning in EMI and CHC environments. In addition to improving academic performance, we would like to improve educator-learner communication, to increase interactive engagement during class time, to develop curiosity with regard to learning new concepts and tools, and to cultivate critical thinking skills.

\subsection{Technology adaptation for instruction} \label{techcas}

In the 21$^\text{st}$-century classroom, embedding an element of technology has become an essential component for successful teaching and learning. Presently, most of our students are `digital natives'. On the other hand, most instructors are `digital immigrants' since their exposure to technology occurred during the latter stages of their lives. Thus, it is hardly surprising that some `digital immigrants' are hesitant in adopting and welcoming technology into their classrooms. `Digital natives' and `digital immigrants' are terms coined by~\cite{prensky2001digital} when he discussed strategies to bridge the communication gap between educators and students. 

\begin{figure}[h]
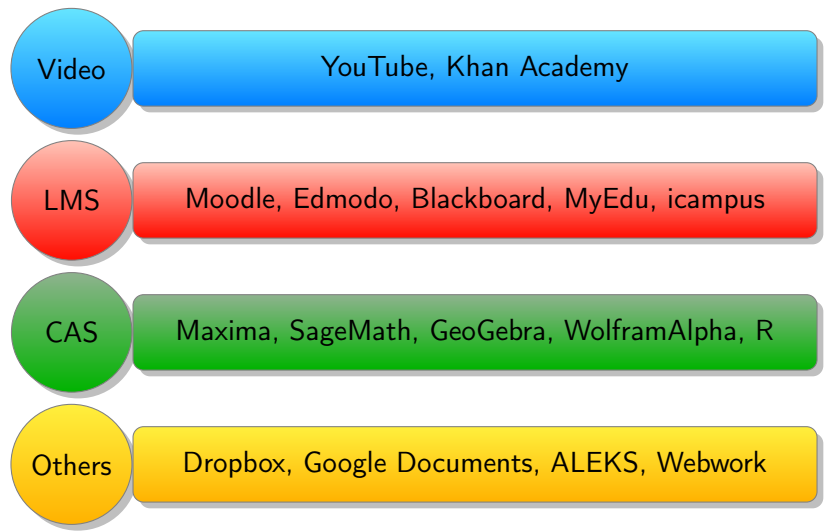

\smartdiagramset{set color list={blue!50!cyan, red!90!orange, green!70!black, red!30!yellow},
  description title font = \sffamily, description font = \sffamily,
  description width = 9cm, description text width = 8cm}
\begin{center}
\smartdiagram[descriptive diagram]{
  {Video,{YouTube, Khan Academy}},
  {LMS,  {Moodle, Edmodo, Blackboard, MyEdu, icampus}},
  {CAS, {Maxima, SageMath, GeoGebra, WolframAlpha, R}},
  {Others, {Dropbox, Google Documents, ALEKS, Webwork}}}
\end{center}
\caption{A few examples of technology that can be utilized for enhancing flipped class activities. {\sl Khan Academy} videos are also accessible via {\sl YouTube}. LMS stands for learning management system and CAS stands for Computer Algebra System. In addition to CAS, technology in mathematics education also comprises numerical computation like {\sl Python}, {\sl Octave}, {\sl Scilab}, interactive between Geometry and Algebra {\sl GeoGebra}, and statistical computing using {\sl R} software.} \label{techno}
\end{figure}

Although we may be labeled `digital immigrants', we are not as hesitant when adopting technology in our classrooms; instead of misgivings, we see opportunities. Our philosophy, instead of forbidding the use of mobile devices during the class, we attempt to utilize and transform the students' attachment to those devices, into a tool for greater engagement while learning mathematics. In general, using technology includes creating our own video recordings, making use of readily accessible {\sl Khan Academy} videos~\citep{noer12}, using an intranet LMS and utilizing technology in mathematics education, for instance using CAS. Some examples of technology that may be incorporated into the flipped classroom are shown in Figure~\ref{techno}.

The use of technology in the classroom has also been explored extensively, but not exhaustively. \cite{davies2013flipping} found that the use of technology in the flipped classroom is effective and facilitated better learning than simulation-based instruction. \cite{kim2014experience} explored flipped classroom blended with technology-enhanced learning, where LMS, {\sl YouTube}, {\sl Google Documents}, {\sl Google Hangout}, {\sl Dropbox} and video cameras are incorporated. See also~\cite{moore2014fostering}. The use of an intelligent tutoring system, known as {\sl Assessment and Learning in Knowledge Spaces (ALEKS)} in a flipped classroom to assist students learning Introductory Statistics has been studied by~\cite{strayer2007effects,strayer2012}. See also~\cite{clark2015effects} where videos and podcasts are incorporated and~\cite{jungic15} where relatively advanced hardware and software are utilized in the video production process. A study of the flipped classroom with the use of technology in mathematics education, using open source mathematical software {\sl GeoGebra} and {\sl Maxima} has been published recently by~\cite{zengin2017investigating}.

Similar to other CAS, {\sl Maxima} is a powerful mathematical software with the capability for the manipulation of symbolic and numerical computations, including differential and integral calculus and matrix operations~\citep{maxima,diehl2009}. {\sl Maxima} is free of charge and is released under the terms of the GNU General Public License. The syntax is relatively easy to grasp by individuals with little or no programming experience. Although another CAS {\sl SageMath} is also popular among university educators because of its user-friendly cloud, the server can be quite lethargic if accessed by individuals with modest Internet connectivity speeds.
For this reason we adopted {\sl Maxima}.

Students in Types A and C were introduced to {\sl Maxima} as part of the initiative to accommodate `digital natives', transform their attachment to technology and to cultivate curiosity while using a new technological tool for improved teaching and learning in SVC. The purpose is not to replace paper-pencil computation proficiency, but rather to assist them, by allowing them to better understand mathematical concepts, develop intuition and to visualize graphical objects, particularly in the context of learning Calculus with flipped classroom pedagogy, as shown by~\cite{heid2001computer},~\cite{karjanto2017adopting} and~\cite{zengin2017investigating}. 

{\sl Maxima} has also been the subject of a number of studies on teaching and learning mathematics using CAS. Here we provide a concise list:
Due to a higher evaluation on portability and aspect of maintainability, {\sl Maxima} is chosen over {\sl Derive}~\citep{garcia2011could}.
\cite{fedriani2011using} present some of {\sl Maxima}'s strengths and weaknesses that they discovered while using {\sl Maxima} to teach mathematics to students majoring in Business. Interactive teaching for Engineering subjects using {\sl Maxima} has been discussed by~\cite{vzakova2011maxima}. Additionally, there are free Calculus e-textbooks incorporating {\sl wxMaxima}~\cite{hannan1,hannan2}.

\section{Research Context: Flipped Classroom in Single Variable Calculus} \label{CA}

\subsection{Course description}

SKKU is a private university affiliated with the Samsung Group and its history dates from 1398. It has two campuses: the Natural Science Campus (NSC) and the Humanities and Social Science Campus (HSCC). SVC, also known as Calculus~1 (course code GEDB001), is offered as a three-credit, General Education, Basic Science and Mathematics (BSM) course. The Department of Mathematics organizes the course as a service course at the NSC and it is also offered at HSCC. Typically, more than 20 sections of the course are offered during the Spring semester while one or two sections are offered in Fall. During Spring, one instructor acts as the course coordinator and communicates to the other instructors on teaching and learning issues, including course topics and the common midterm and final examinations. During Fall, instructors organize and administer their own assessments and there are no common examinations.

SVC is designated as an EMI class and the primary target audience is freshmen, both local and international students. Generally, the content in SVC in a South Korean university is equivalent to the content covered in both Calculus~1 and Calculus~2 in an American university. On the other hand, Calculus~2 (SKKU course code GEDB002) at SKKU and other tertiary institutions in South Korea is typically equivalent to the Calculus~3 course offered in many American universities. For both courses, the adopted textbook is~\cite{stewart2010}. Many translations are available, including Korean; students have been observed to and may freely use any translation to study from in and out of class. Nevertheless, all teaching materials and student assessments are provided in English. Students submit and write their assessments, including their examinations, in English. We emphasize having a good command of English in all of our teaching and learning activities, both oral and written form, the latter is communicated via the intranet LMS {\sl icampus}.
\begin{table}[h]
\begin{center}
{\small
\caption{\small Various types of instruction and course logistics for flipped classroom investigated in this study.} \label{fliptype}
\begin{tabular}{@{}lllll@{}}
\toprule
			            & \multicolumn{4}{c}{Instruction type} \\ \cline{2-5}
			            & \qquad A  			   & \qquad B 			  & \qquad \quad C  	    & \qquad D \\ \hline
Description  			& fully flipped 		   & fully flipped 		   	  & single-topic flipped & unflipped \\ \hline
Video source 			& instructor-created       & third-party              & instructor-created   & no video \\ \hline 
Number of videos    	& 12                       &  10                      &  2                   &  0 \\ \hline 
Technology used     	& {\sl YouTube}, LMS, CAS  & {\sl Khan Academy}, LMS  & {\sl YouTube}, LMS, CAS   & LMS \\ \hline 
\multirow{2}{*}{Pedagogy}& in-class problem  	   & in-class problem  		  & mixed lecture and    & traditional \\
                         & solving         	 	   & solving                  & problem solving      & lecture \\ \hline
\multirow{2}{*}{In-class activities} 	
						& quiz and                 & small group              & note taking and   	 & note taking \\
				        & problem solving          & problem solving          & problem solving      &  \\ \hline
\multirow{2}{*}{Out-of-class activities} 
						& online video and 	       & online video and		  & online video and     & homework \\
					    & reading assignment       & reading assignment		  & group project	     &	\\ \hline
Type of exam            & common exam              & common exam              & individual exam      & common exam \\ \hline
Semester		        & Spring 2016              & Spring 2015              & Fall 2015            & Spring 2015 \\ \hline
Instructor (Section)    & NK (57)                  & LS (50)                  & NK (42)              & LS (49) \\
\bottomrule  
\end{tabular}
}
\end{center}
\end{table}

Flipped classrooms in this study were conducted over a three-semester period from Spring 2015 to Spring 2016. During Spring 2015, author (LS) taught two sections of SVC, one section is conducted with flipped classroom instruction (Type B) and the other is delivered in the traditional lecture format (Type D), which acts as a control. While NK taught one section of SVC each during Fall 2015 and Spring 2016, NK adopted flipped classroom pedagogy for the entire Spring 2016 semester (Type A) and flipped the topic Introduction to Differential Equations in Fall 2015 (Type C). For the rest of the paper, these four sections are indicated by instruction Types A, B, C and D.  Both Types A and B are flipped classrooms where the pedagogy is conducted throughout the semester; differences lie in the source of the video recordings and the use of technology. Type A's videos are instructor-created while Type B's videos are third-party created, i.e. {\sl Khan Academy}. NK also introduced {\sl Maxima} in Types A and C throughout the semester. The four types of instruction are summarized in Table~\ref{fliptype}.

\subsection{Technology utilization in flipping}

As mentioned in Subsection~\ref{techcas} of the study's theoretical framework, an adoption of technology in teaching and learning is an essential component for implementing flipped pedagogy successfully. In our flipped classrooms, the adaptation of technology includes video recordings, using an intranet LMS {\sl icampus} and the {\sl Maxima} for Types A and C.
\begin{figure}[h]
\begin{center}
\includegraphics[width = 0.45\textwidth]{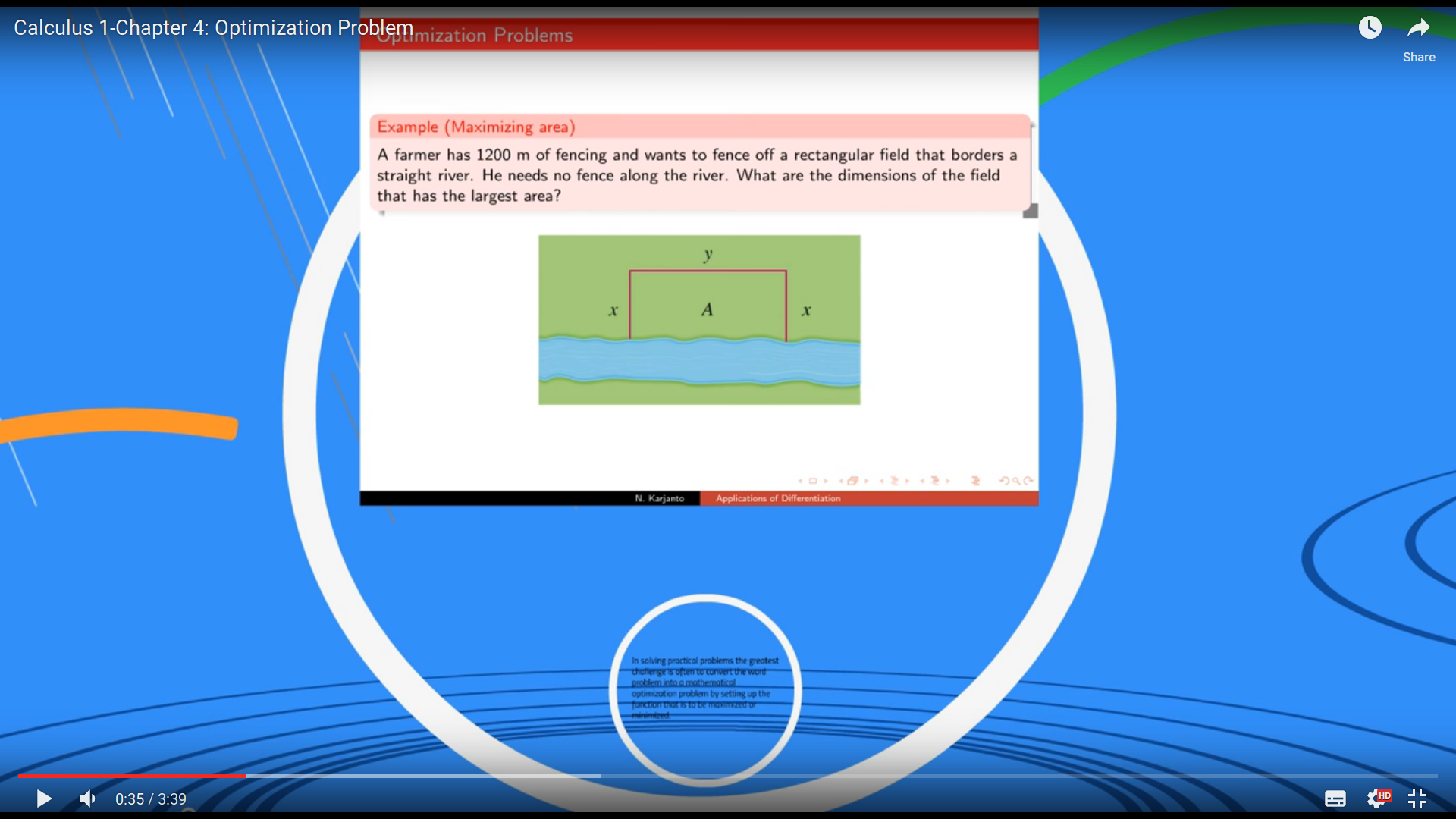} \hspace{1cm}
\includegraphics[width = 0.45\textwidth]{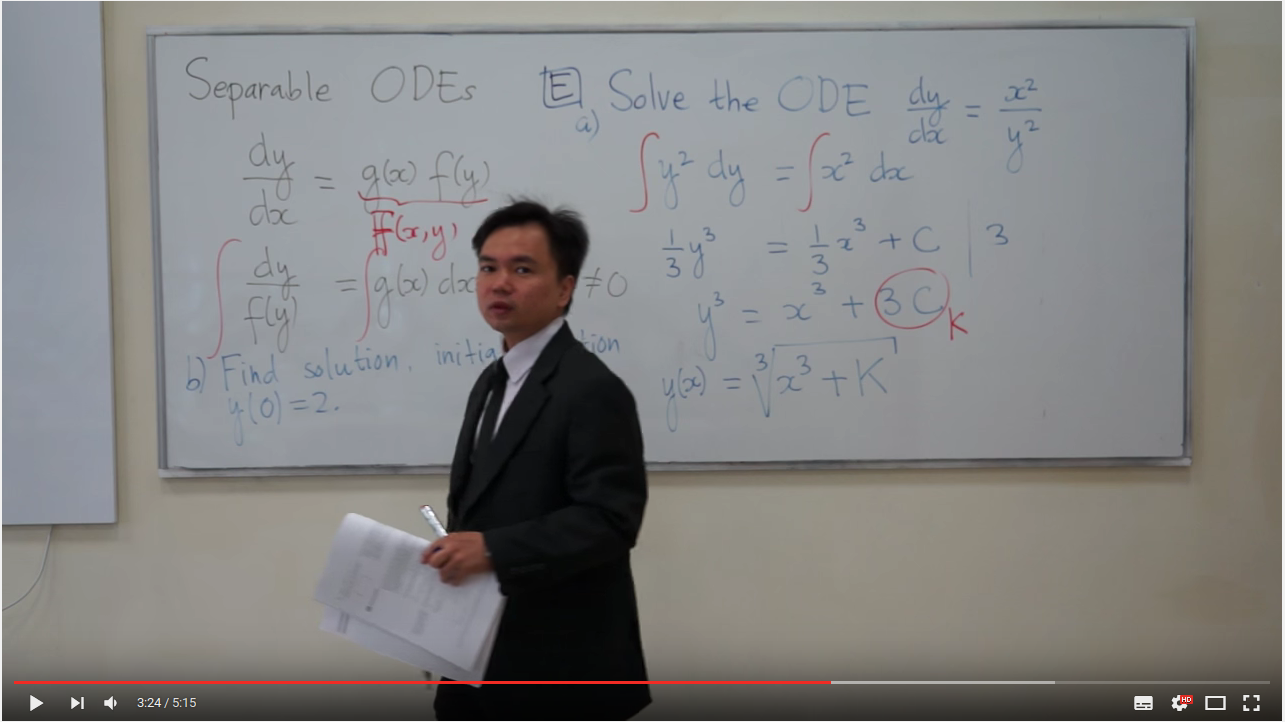} \vspace{0.5cm}\\
\includegraphics[width = 0.45\textwidth]{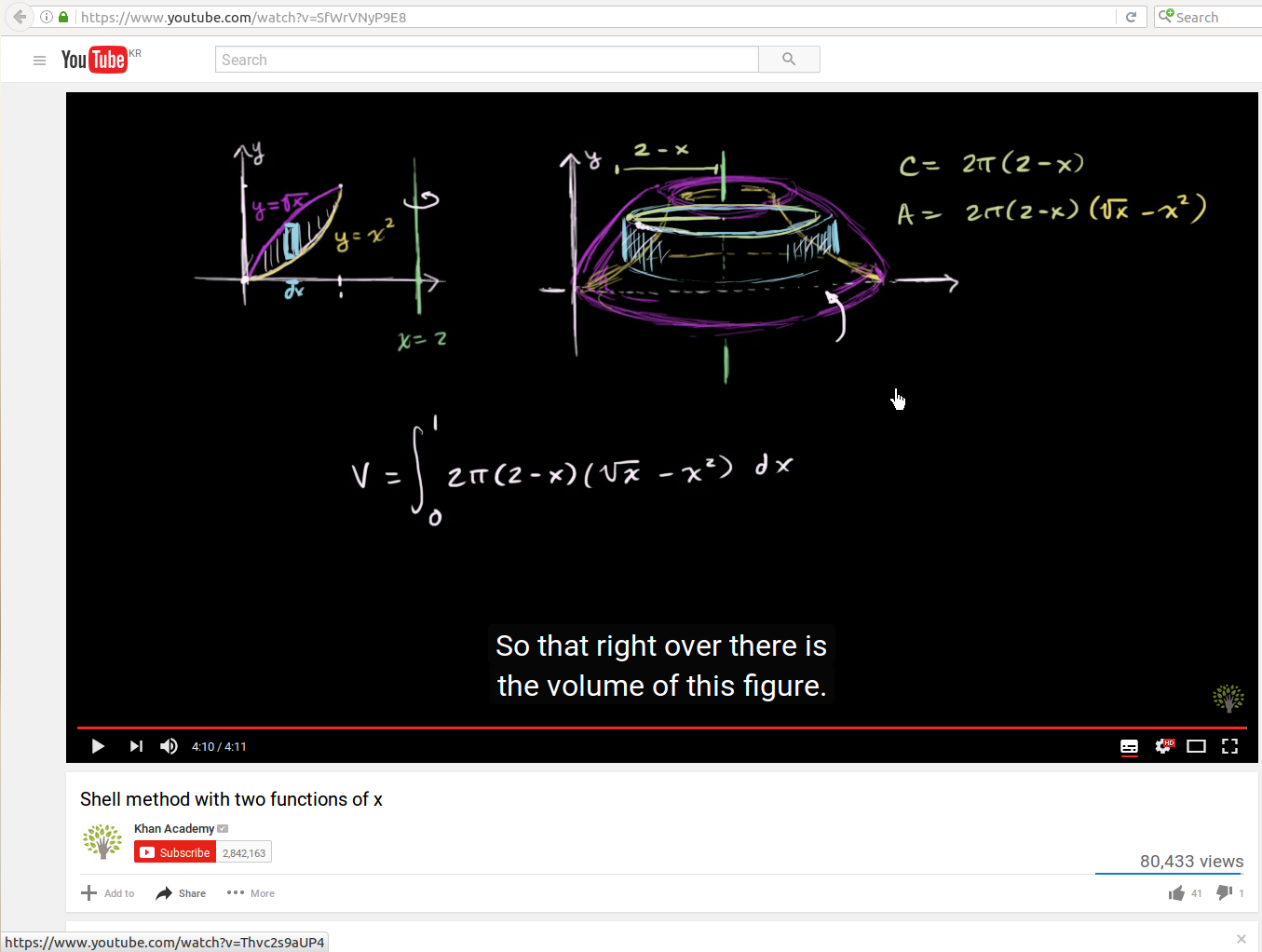} \hspace{1cm}
\includegraphics[width = 0.45\textwidth]{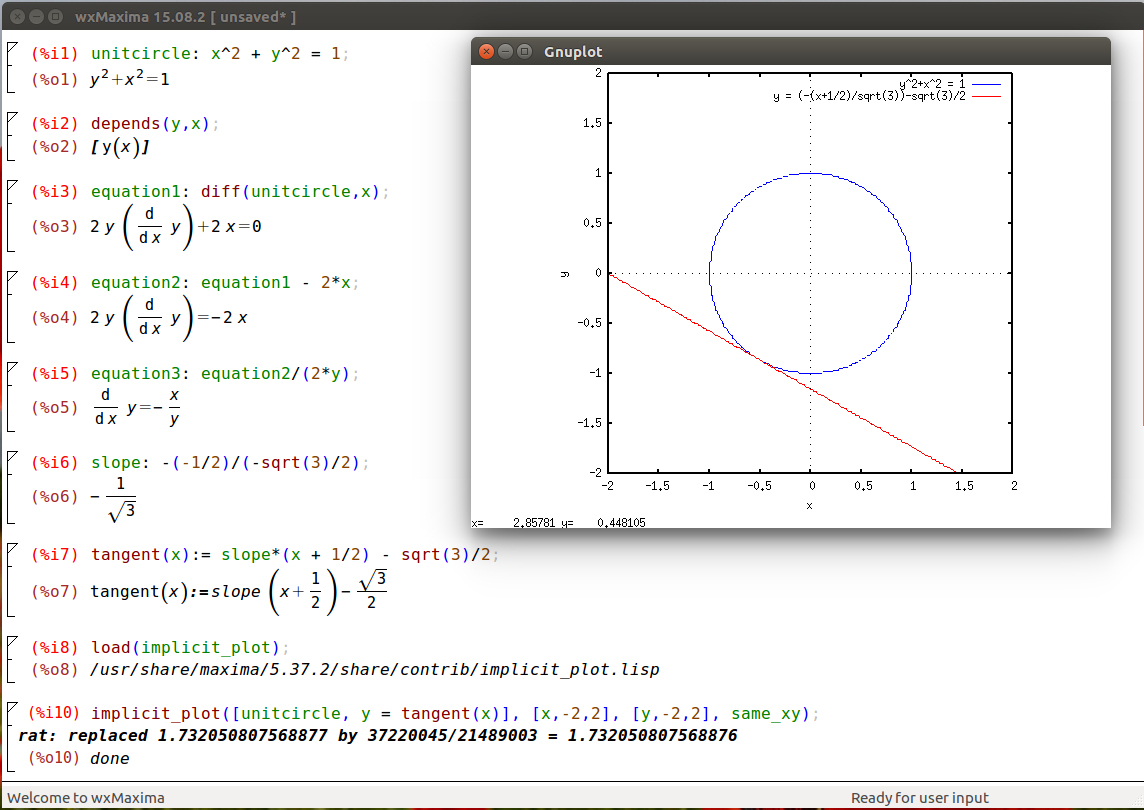} \\
\end{center}
\caption{\small (Top panels) Samples of recorded video lectures used in Types A and C. The left screenshot explains an optimization problem for maximizing the area of a rectangular field by fencing it. The right screenshot deals with the method of integrating factor in Differential Equations. The durations of these videos are 3:39 and 5:15 minutes, respectively. (Bottom left) A sample of video from {\sl Khan Academy} with duration of 4:11 discussing the calculation of the volume of a solid revolution using the method of cylindrical shells. (Bottom right) A~screenshot of {\sl Maxima} of a problem in finding a tangent line of a unit circle at the point $(-1/2,-\sqrt{3}/2)$.}
\label{videos}
\end{figure}

All recorded video lectures were made accessible via {\sl YouTube} for Types A and C. For Type B, {\sl Khan Academy} video recordings are also accessible via {\sl YouTube}. Depending on the material, the duration of the recorded videos ranged from 10~to~15~minutes. An empirical study on massive online open courseware suggests that videos with a shorter duration are more engaging to students than those of longer duration~\citep{guo2014video}. Some of Type A's videos were produced using {\sl RecordMyDesktop} software. There is no camera to capture the face of the instructor speaking and explaining the material, so Type A's videos are similar to {\sl Khan Academy} videos, albeit the latter are professionally prepared. Additionally, the two recordings included in Type C were captured using a digital single-lens reflex camera. The instructor's face and the upper body are visible in both videos. Screenshots examples of the screencast video recording are shown in Figure~\ref{videos} and URL links are sent to the students using the LMS {\sl icampus}. The presentations were prepared using {\sl Beamer} and {\sl Prezi}.

For Type A, some in-class activities include demonstrations of simple calculations and curve sketching using {\sl Maxima}. The students are also encouraged to try the CAS as part of problem-solving activities, by checking and confirming answers obtained by hand calculation using the CAS and to visualize graphs of function. One of the out-of-class activities using {\sl Maxima} involved finding and sketching a tangent line at one point of an implicit function, in this case, a unit circle. A video demonstrating a similar problem using a folium of Descartes had been recorded and the URL link was sent to the students to view via the LMS. The {\sl Maxima} source code and resulting sketch of the corresponding assignment are displayed in the lower right panel of Figure~\ref{videos}. For Type C, {\sl Maxima} is only utilized for demonstration during the lecture and zero CAS related assignments were given.

\section{Methodology} \label{method}

In this study, we adopted an observational research method; where as instructors, we observed the effect of implementing the flipped classroom pedagogy on our students. A case study is applied to four sections of the SVC course with different instruction types. We also conducted a mixed methods study by applying both quantitative and qualitative methods.

For the quantitative aspect, we compared students' academic achievement. Data for academic performance are obtained from midterm and final examination scores, as well as from the end of semester (final) letter grade distribution. The latter includes homework, assignment (project), and quiz assessments. In particular, we are interested in discovering whether different types of flipped instruction reveal a statistically significant difference in terms of academic performance. The qualitative method is used to share students' perception of their flipped classroom experience.

\subsection{Participants}

The participants for this study are students enrolled in the four sections of our SVC classes offered at SKKU's NSC, from Spring 2015 to Spring 2016. Three classes were designated as experimental groups (Types A--C) and another class acted as the control group (Type D). The dataset consists of 310 students, 225 were distributed among the three experimental groups while the remaining 85 belonged to the control group. Participants' ages ranged from 18 to 24 years old. The convenience sampling method was employed in selecting the study groups due to accessibility and efficiency. Detailed demographics are displayed in Table~\ref{demog}.

\begin{table}[h]
\begin{center}
\caption{\small The demographics of the participants. International origin refers to both Korean students who graduated from secondary schools abroad and non-Koreans irrespective of their high school origin.} \label{demog}
\begin{tabular}{@{}llrrrrr@{}}
\toprule
            &               	& \multicolumn{4}{c}{Instruction type} & \\ \cline{3-6}
            &               	& {\qquad} A & {\qquad} B\, & {\qquad} C\, & {\qquad} D\, & \qquad Total  \\ \hline
Enrolled    &               	& 59 & 85 & 81 & 85 & 310 \\ 
\multirow{2}{*}{Gender}
            & Male          	& 50 & 65 & 48 & 67 & 230 \\
            & Female        	&  9 & 20 & 33 & 18 &  80 \\ \hline
\multirow{5}{*}{Class Year {\quad}}  
			& Freshman 		    & 48 & 78 & 60 & 78 & 264 \\
            & Sophomore {\quad} &  9 &  3 & 10 &  3 &  25 \\
            & Junior            &  0 &  2 &  4 &  2 &   8 \\ 
            & Senior        	&  0 &  0 &  0 &  0 &   0 \\
            & Post-senior       &  2 &  2 &  7 &  2 &  13 \\ \hline
\multirow{2}{*}{Status}      
            & First-time    	& 55 & 85 & 53 & 85 & 278 \\  
            & Retake        	&  4 &  0 & 28 &  0 &  32 \\ \hline
\multirow{2}{*}{Origin}      
			& International 	&  6 &  5 & 31 &  4 &  46 \\
            & Local         	& 53 & 80 & 48 & 81 & 262 \\            
\bottomrule
\end{tabular}
\end{center}
\end{table}

\subsection{Measurement}

Students' academic achievement is measured by two major exams and final letter grades. The former consists of one-hour long midterm and final examinations. The exam questions are presented in Appendices A--C. Each examination's maximum possible score is~100. The final grade is comprised of 80\% from the two examination and 20\% from other assessments; the letter grades are distributed according to a relative grading policy~\citep{handbook}.

The four class types distributed the 20\% among homework, quizzes and other assignments. Type A's distribution is quizzes, 10\% and online {\sl MapleTA} assignments, 10\%. Type B's distribution includes two quizzes, three homework sets and class attendance. On the other hand, for Type C, it consists of a project assignment with a written report on an application of Calculus in Economics. For Type D, it is composed of two quizzes, four homework sets (paper-based and {\sl MapleTA}) and class attendance. Furthermore, the letter grades range from the lowest, F for failure, to the highest, A$+$. The letter grade and its equivalent value grade point average (GPA) calculations are A$+ \equiv 4.5$, A $\equiv 4.0$, B$+ \equiv 3.5$, B $\equiv 3.0$, C$+ \equiv 2.5$, C $\equiv 2.0$, D$+ \equiv 1.5$, D $\equiv 1.0$ and F $\equiv 0$.

Students' perception of the flipped classroom was collected through a set of online questionnaires administered by the Registrar Office. These are accessible through an intranet university system for faculty and students, known as {\sl Advanced Sungkyunkwan Information Square--Gold Lawn Square} (ASIS-GLS). 
The online questionnaires were administered twice during the semester; on both occasions one week before the examination periods. The first questionnaire contains two open-ended questions related to teaching and learning: ``Please write down the most impressive and instructive aspects of this class.'' and ``Please write down your suggestions for the instructor to improve the class.'' The second questionnaire consists of six statements on a five-item Likert scale and a single open-ended question, identical to the second question from the first questionnaire. The the six statements are part of a standard evaluation for general teaching and are not specifically related to flipped classroom pedagogy, thus we will not discuss these results here. Rather, we will focus on the students' perceptions about the flipped classroom gathered from their feedback. 

\subsection{Data collection and analysis}

Students' demographic data are collected from the Registrar Office. A~number of simple statistical analyses were carried out in this study including the Student's $t$-test, one-way analysis of variance (ANOVA) and Tukey's honestly significant difference (HSD) test. Data analysis was conducted using {\sl LibreOffice Calc} and {\sl GeoGebra}. ANOVA was conducted online, a statistical computation {\sl VassarStats} along with Dr. Lawrence Turner's website hosted by Southwestern Adventist University, Keene, Texas. The box plots and histograms were generated by R~software~\citep{rproject15}.

\section{Results} \label{finding}

This section discusses quantitative results related to academic performance and qualitative results sourced from our SVC students' perceptions of flipped classroom pedagogy.

\subsection{Academic performance} \label{quanti}

\begin{figure}[h]
\begin{center}
\includegraphics[width = 0.45\textwidth]{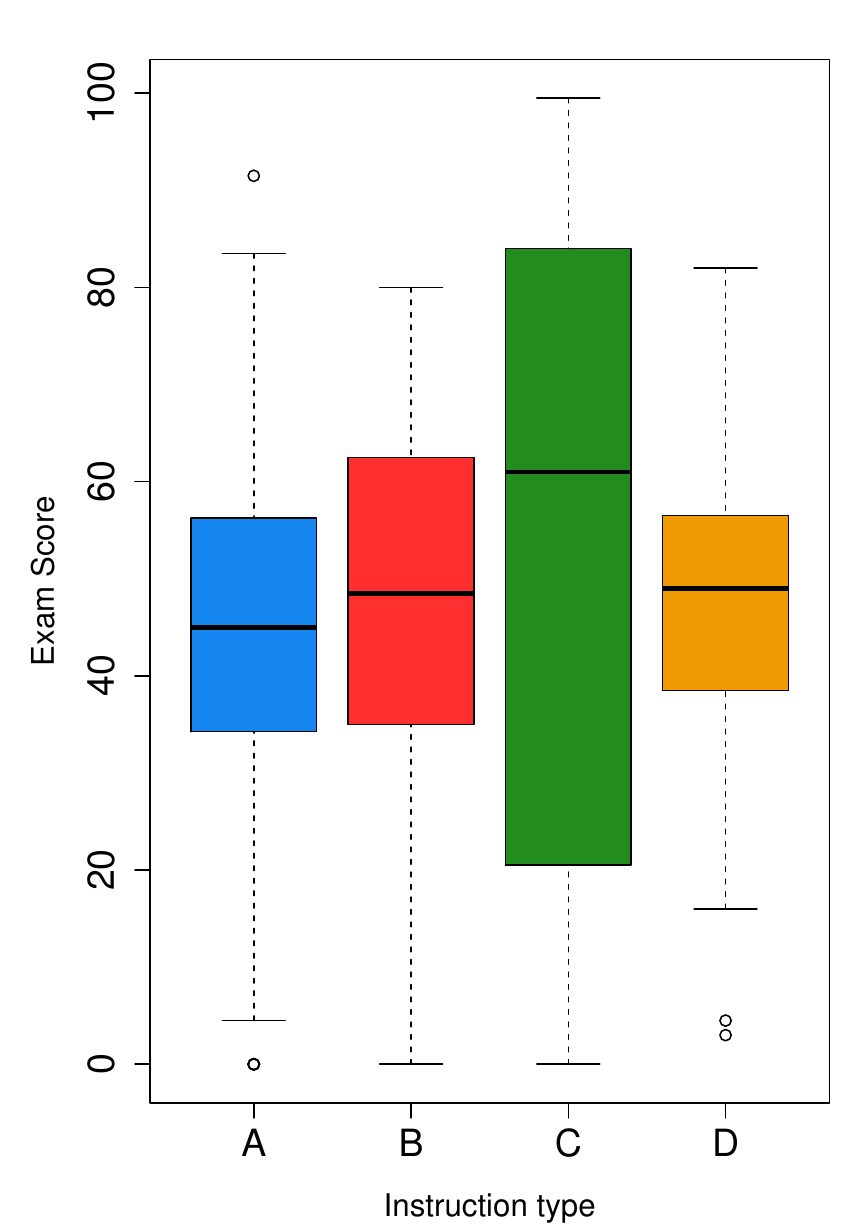}    \hspace{1cm}
\includegraphics[width = 0.45\textwidth]{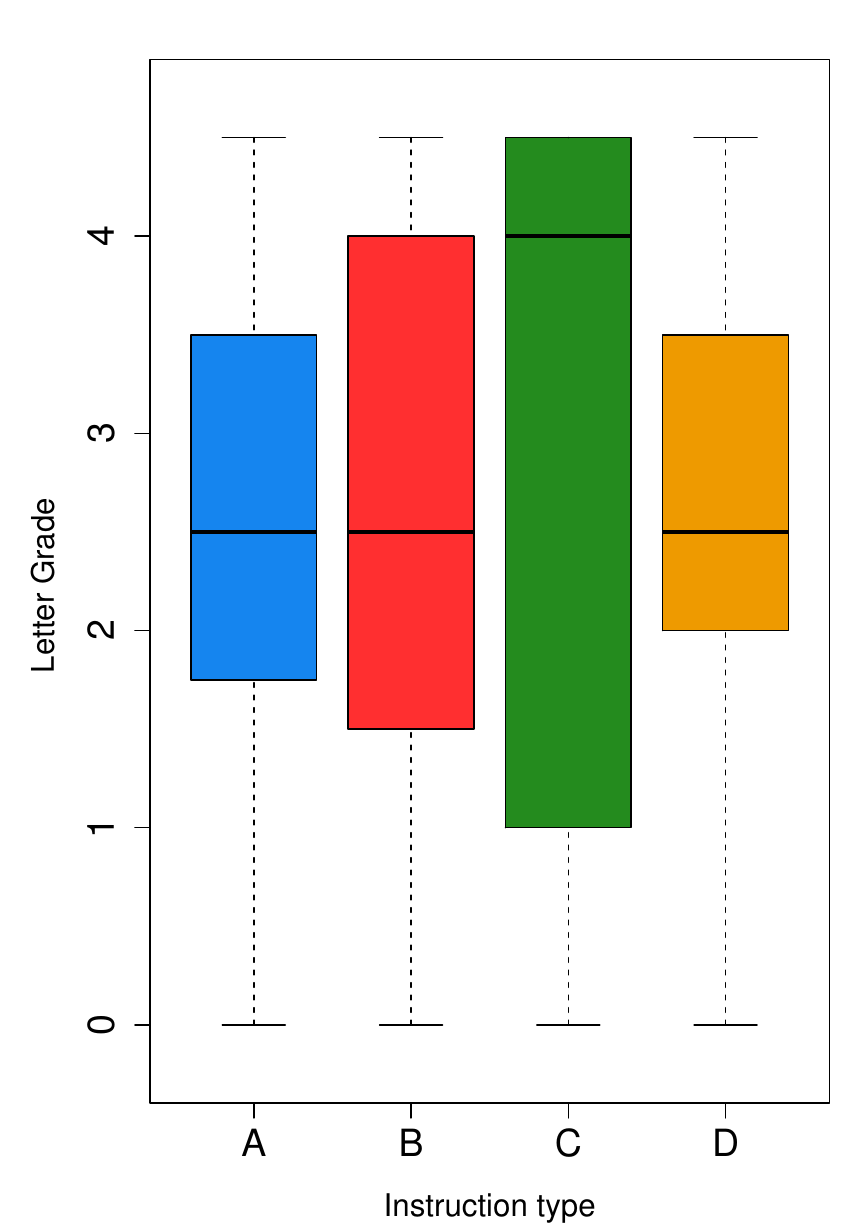} 
\end{center}
\caption{ANOVA box plots of the exam score (left panel) and the final letter grade distribution (right panel)
for flipped classroom Type A (blue box), Type B (red box), Type C (green box) and traditional classroom Type D (orange box).} \label{boxplot}
\end{figure}

\begin{figure}[h]
\begin{center}
\includegraphics[width = 0.9\textwidth]{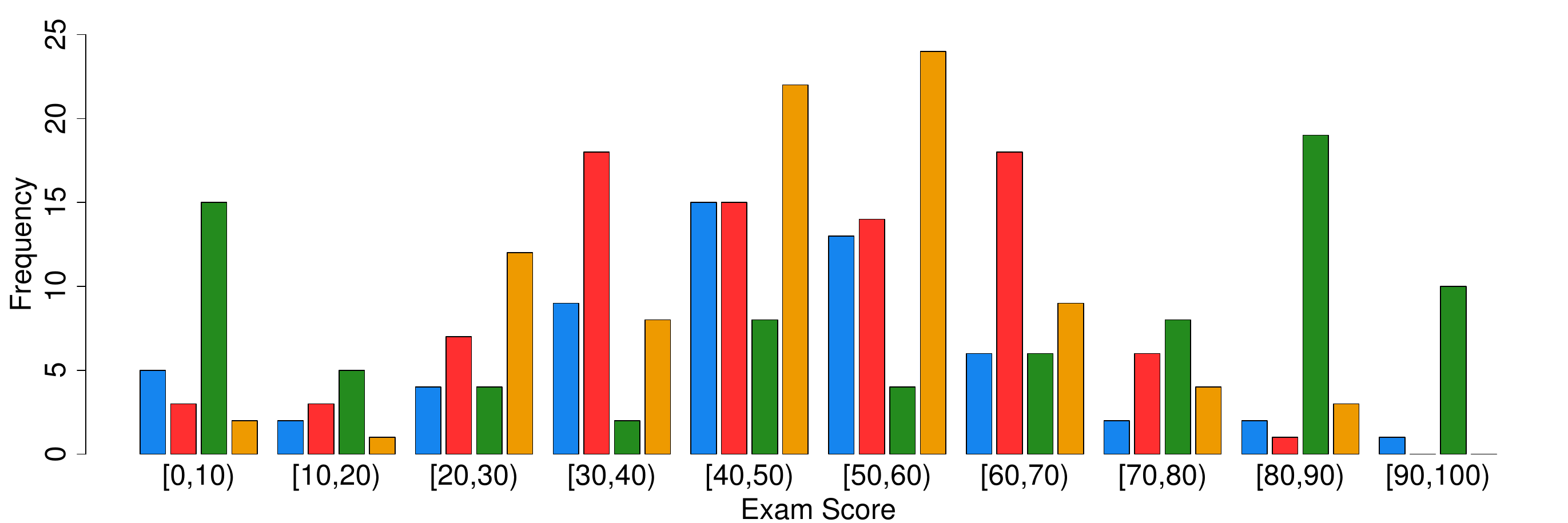}    \\
\includegraphics[width = 0.9\textwidth]{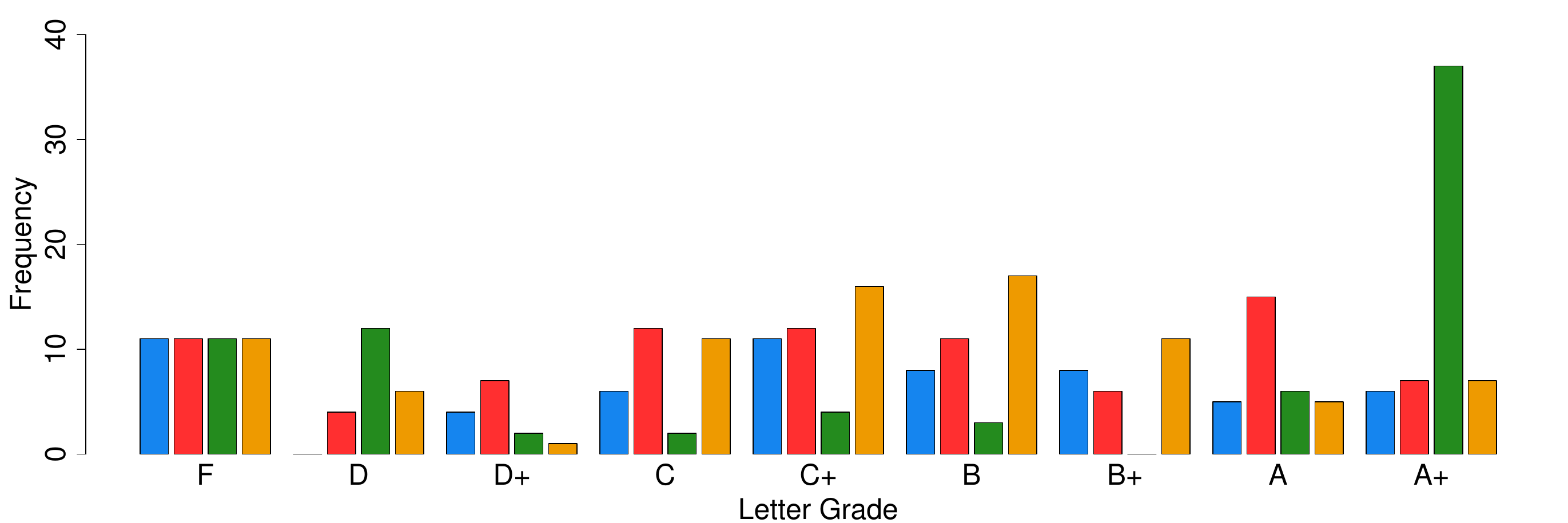} 
\end{center}
\caption{Bar plots of the exam score (top panel) and the final letter grade distribution (bottom panel)
for flipped classroom Type A {\color{black} (blue bars)}, Type B {\color{black} (red bars)}, Type C {\color{black} (green bars)} and traditional classroom Type D {\color{black} (yellow bars)}.} \label{barplot}
\end{figure}

This subsection explains quantitative findings of students' academic performance based on assessment results. 
The assessment data include the examination score, which is the average of midterm and final scores and the end of the semester final grade which is in the form of letter grades. Means and standard deviations of assessment results, both for the exam score and the letter grade distribution are displayed in Table~\ref{scoregpa}. Furthermore, box plots and histograms for both the exam score and letter grade distributions are shown in Figures~\ref{boxplot}~and~\ref{barplot}.
\begin{table}[h]
\begin{center}
\caption{\small Means and standard deviations of the assessment result for various types of instruction and the number of students in each class. Data for examination score (average of the midterm and final exam scores, out of 100) and final letter grade (out of 4.5) are displayed in Columns 3--4 and 5--6, respectively.} \label{scoregpa}
\begin{tabular}{@{}ccccccc@{}}
\toprule 
\multirow{2}{*}{Instruction type} & \multirow{2}{*}{$\qquad N \qquad$} & \multicolumn{2}{c}{Examination score} & & \multicolumn{2}{c}{Letter grade} \\ \cline{3-4} \cline{6-7}  
            &     & Mean	   & SD  			 &    & Mean         & SD      \\ \hline
A           & 59  & 44.170     & 20.621          &    & 2.449        & 1.435   \\ 
B           & 85  & 46.950     & 18.106          &    & 2.518        & 1.372   \\ 
C           & 81  & 54.451     & 33.415          &    & 2.969        & 1.745   \\ 
D           & 85  & 47.547     & 15.939          &    & 2.477        & 1.293   \\ \cline{1-2} 
Total       & 310 &            &                 &    &              &         \\
\bottomrule  
\end{tabular}
\end{center}
\end{table}

We applied ANOVA analysis to both the exam score and to the final letter grade under an assumption of independent samples with a standard weighted-means analysis. The ANOVA summary of academic performance across the four types of instruction is displayed in Table~\ref{examlett}. For the exam score means, there was a statistically significant difference across the four types of instruction as determined by the one-way ANOVA ($F(3,306) = 2.67$, $p = 0.0477$). Further, the effect size value ($\eta^2 = 0.0255$) suggested a small practical significance. By considering that our data meet the variance homogeneity assumption, we adopted Tukey's honestly significant difference (HSD) post-hoc test. The post-Tukey HSD test of the ANOVA analysis for the exam score indicates that there is a significant difference between the exam score means of Types A and C, with $p$-value $= 0.0477 < 0.05$. There are no significant differences in the exam score means between instruction Types A and B, Types A and D, Types B and C, Types B and D, as well as between Types C and D. Furthermore, the absolute (unsigned) difference between any two sample means required for significance at the designated level is given as follows: HSD$(0.05) = 9.7$ and HSD$(0.01) = 11.78$. On the other hand, the ANOVA analysis of the final grade indicates that there is no significant difference in terms of the letter grade means across all types of instruction as determined by one-way ANOVA ($F(3,306) = 2.19$, $p = 0.0892$). Further, the effect size value ($\eta^2 = 0.0210$) suggested a small practical significance.
\begin{table}[h]
\begin{center}
\caption{One-way ANOVA summary of academic performance result, both exam score and final letter grade.}	 \label{examlett}
\begin{tabular}{@{}llrrrccc@{}}
\toprule
          	 & Source 		  & SS         & df  & MS      & $F$  & $\eta^2$ & $p$    \\ \hline 
Exam score   & Between groups &    4255.43 &   3 & 1418.48 & 2.67 & 0.0255   & 0.0477 \\
             & Within  groups & 162,864.27 & 306 &  532.24 &      &          &        \\ 
             & Total          & 167,119.70 & 309 &         &      &          &        \\ \hline  
Letter grade & Between groups &      14.21 &   3 &    4.74 & 2.19 & 0.0210   & 0.0892 \\
             & Within  groups &     661.70 & 306 &    2.16 &      &          &        \\
             & Total          &     675.91 & 309 &         &      &          &        \\
\bottomrule
\end{tabular}
\end{center}	
\end{table}

Two possible explanations for the statistically significant difference in the exam score between Types A and C exist. First, there is a significant number of repeaters enrolled in Type C, 28 out of 81, which accounts for almost 35\% of the class population. On the other hand, the number of repeaters enrolled in Type A is only four and there are zero repeaters in Types B and D. Repeaters might have an advantage over their first-time peers because of past exposure to and study of the material. This may translate to a higher average exam score among the retake group and hence a higher average exam score for Type C. Second, the level of difficulty of the exam questions given to Types A and C is slightly different. As mentioned in Table~\ref{fliptype}, all instruction types, except Type C, administered common exams for both the midterm and final. The common exams, generally have seven to nine questions, depending on whether the questions require longer or shorter solutions. The questions and format are new for many students since the questions are neither given in homework assignments and quizzes nor taken directly from the textbook. In contrast, the exam administered for Type C contained five questions each in the midterm and final exams. In addition, some of these questions were discussed previously during in-class problem-solving sessions.

On the other hand, there is no significant difference for the final letter grades across the four types of instruction. The best explanation for this may be that a relative grading rule is applied when assigning the letter grades. See Figure~\ref{barplot}, the bottom panel and notice that Type C has many high achievers. Interestingly, it was discovered after the midterm that this particular class was composed of students with a diverse range of mathematical skills and abilities. Generally, we observed that the students who received their secondary education in South Korea had little to no difficulty following the course. On the other hand, the students who received their secondary education outside of South Korea, irrespective of their nationalities, had some difficulties following the course. 

A $t$-test result for difference of means in academic performance for Type C, excluding repeaters, shows a very significant difference between non-repeaters and repeaters. Out of 53 students (81 total students subtract 28 repeaters), 25 students graduated from South Korean secondary schools, while 28 students graduated from various international secondary schools. Note that the 31 students with international origin displayed in Table~\ref{demog} include these 28 students and international students who graduated from South Korean secondary schools. Table~\ref{ttestlocint} displays a summary of the $t$-test results for both exam score and the final letter grade distribution average. In both instances, the $p$-value $= 0 < 0.01$ (extremely significant) and Cohen's effect size values ($d = 1.5564$ and $d = 1.3242$) suggested very large practical significances for both exam score and letter grade, respectively.
\begin{table}[h]
{\small
\caption{Summary of $t$-test for difference of means in academic performance between students attended secondary school locally or internationally. The exam score and the letter grade distribution are considered and in both instances, the effect sizes are very large.} \label{ttestlocint}
\begin{center}
\begin{tabular}{llcrrccccc}
\toprule
					        & High school   & N  & $\bar{x}\; \quad$ & sd {\quad}      & SE     & df      & $t$    & $p$ & Cohen's $d$ \\
           \cline{2-10}
\multirow{2}{*}{Exam score} & Local         & 25 & 68.3200   & 26.4964 & \multirow{2}{*}{7.4357} & \multirow{2}{*}{50.7163} & \multirow{2}{*}{5.6627} &     \multirow{2}{*}{0} & \multirow{2}{*}{1.5564}             \\
				            & International & 28 & 26.2143   & 27.6005 &  &  &  & &       \\
\hline         

\multirow{2}{*}{Letter grade} & Local         & 25 & 3.6400   & 1.4686 & \multirow{2}{*}{0.4179} & \multirow{2}{*}{50.8884} & \multirow{2}{*}{4.8126} &     \multirow{2}{*}{0} & \multirow{2}{*}{1.3242}             \\
				            & International & 28 & 1.6250   & 1.5731 &  &  &  & &       \\
\bottomrule  
\end{tabular}
\end{center}
}
\end{table}

\subsection{Students' perception and feedback} \label{quali}

This subsection discusses students' perceptions and feedback after implementing flipped classroom pedagogy. The feedback is arranged in chronological order; first, Type B (Spring 2015), followed by Type C (Fall 2015) and then Type A (Spring 2016). Type D is not delivered as flipped classroom pedagogy, thus the corresponding feedback is excluded.

\subsubsection{Students' feedback--Spring 2015}

In-class activity for Type B is typically structured as follows: after a brief introduction to begin each class, which includes an overview of the material, problems related to the overview are distributed and attempted. The problems are selected from the course textbook, other Calculus textbooks and past examination papers. The students are given time to work individually or cooperatively in teams; the teams are assigned but the team assignments are not strictly enforced. Students may also seek assistance from the instructor during this time period. Afterwards, students are invited to vote for those questions they need explained; the questions are ranked by popularity (need) and the presenter (either students or the instructor) presents solutions.

Survey results via {\sl Google Forms} indicate that 40 out of 45 students favored this class structure, while 5 preferred the structure without an overview. Regarding assigned {\sl Khan Academy} videos for viewing, two students viewed often, 12 students watched sometimes and 31 students viewed infrequently or did not view at all. Similar results were found for the reading assignments, three students read often, 19 students read sometimes and 23 students responded seldom or not at all. In the context of EMI, no feedback was obtained in connection with the ease or difficulty involved in following classes in English. Nonetheless, 60 students wrote their comments in Korean, it seems that around 70\% of the class population is more comfortable providing feedback in their native language. 

Below are some positive comments from the students in Type~B:
\vspace*{-0.3cm}
\begin{quoting}[leftmargin=0.0\parindent]
{\small ``I think the flipped classroom system was very efficient because I was motivated to self-study the materials before coming to the lectures, and most importantly, there was no homework. Due to that, I could have a self-driven pace of learning, and revise what I really needed to.'' \par}  \vspace*{0.1cm} 
{\small \setlength{\parindent}{0cm} ``This class is different from other classes because the professor and the students communicate a lot during the class. I think it's a good idea! Professor LS listens to us very carefully. I love it! Professor LS is so kind to us.'' \par}
\end{quoting}
\vspace*{-0.3cm}

As we can see from these comments, one student praised the efficiency of Type~B because he/she felt more motivated to do independent study before coming to the classroom. The student also remarked that there was no homework (in a traditional sense); but did not realize that the homework being completed in class. Another student noticed that more improved and communication exists between the instructor and the students in this class in contrast to other classes.

The students in Type~B made the following suggestions:
\vspace*{-0.3cm}
\begin{quoting}[leftmargin=0.0\parindent]
{\small ``Although students were grouped, there were not much group activities. I think if there is some element of
competition among groups with some form of a rewarding system, it will improve participation of students during lecture.'' \par} 
{\small \setlength{\parindent}{0cm} ``It would be great if she could teach the basics first, instead of assuming that most of us know.''}
\end{quoting}
\vspace*{-0.3cm}

One student pinpointed a defect in the team activities; although the purpose of the teams for small group activities is to encourage collaboration among the students, this student felt that competition among peers should be kept. He/she also suggested that if there is a reward system for team activities, there would be increased activity among team members and among teams in the class. Another student preferred a type of mini-lecture at the beginning of the class to discuss the basic materials before plunging into problem-solving activities. This is a pertinent suggestion, but it would be interesting to identify this student's survey response regarding video viewing, often, sometimes or seldom, not at all; since the basic materials, definitions, theorems and worked examples are explained thoroughly in the videos. 

\subsubsection{Students' feedback--Fall 2015}
The SVC course during Fall 2015 was not officially offered as a flipped classroom pedagogy course. In both midterm and final semester feedback, the common comment thread was about there being too much material to cover during a limited time period. Many felt that the instructor covered the material too quickly and not thoroughly. They preferred slower step by step explanations when discussing particular examples or problems from the textbook. One student suggested that it would be great if the course could be offered as a flipped classroom: {``I think you will be a great professor if you open flipped classes.''}

Since it required more preparation if we tried to implement flipped classroom pedagogy during the second half of the semester, instead of full implementation, we decided to attempt flipping only one topic, Introduction to Differential Equations, i.e. Chapter 7~\citep{stewart2010}. Unfortunately, there is no feedback related to this minor attempt of flipping at the end of the semester, so we do not know whether the students enjoy a flipped approach on either individual or even this specific topic. 

In the context of EMI, there is no comment in relation to the difficulty of following classes in English. We also observed that students of international origin have no difficulty understanding English. However, 29 and 33 students responded to open-ended inquiries in Korean during the midterm and the final course evaluation, respectively. It seems 40\% of the student population in Type C are more comfortable expressing their thoughts in Korean when providing feedback. No feedback was given in relation to adopting {\sl Maxima} during teaching and learning sessions. This is most likely due to CAS being used for demonstrations and not as part of an assessment.

\subsubsection{Students' feedback--Spring 2016}
During Spring 2016, only one out of 22 sections SVC offered in NSC was designated as a flipped classroom, referred as Type A in this paper.
Selected students' comments on the pedagogy during mid-semester feedback are given as follows: the first two are positive and the last two are negative.
\begin{itemize}[leftmargin=1em]
{\small
\item Conducting problem-solving in advance and preparing by means of a flipped class are a good way to study.
\item The instructor has made an effort to conduct flipped class very well in a progressive manner.
\item Learning on video and practicing in the class doesn't make sense.
\item I would like to receive more detailed instruction on the video recordings. 
It's okay if the duration will be longer than the current video recordings. 
I have a trouble with understanding something. \par}
\end{itemize}

We observe that the students' reaction is mixed. Positive comments include class preparation prior to attending classes and an appreciation for the instructor's effort to implement the pedagogy progressively. On the other hand, one student felt that this approach does not really assist him/her in his/her learning. Another student seems to encounter difficulty in understanding the materials presented in the video. This might be caused by lower quality instructor-created videos as opposed to high quality {\sl Khan Academy} videos, or poor student English or mathematical comprehension skills.

In the context of EMI, 36 and 18 students wrote their comments in Korean during the midterm and the final periods, respectively. This shows that up to 60\% of the student population in Type A are more comfortable writing in Korean instead of English. Although there is no feedback related to difficulty in following classes in English, we often encountered students who inquired whether they can write their assessments in Korean instead of in English. In connection to incorporating {\sl Maxima} while teaching and learning, some students expressed difficulty in using it because the software is not familiar to them. They also would like to obtain further explanation not only in operating it but also in the technical aspects of downloading and installation.

In particular, regarding one assignment related to the CAS, one student sent an electronic message to indicate that he/she encountered difficulty while attempting to complete the assignment. The email was sent on 20 March 2016. Name and student ID are withheld for privacy. Some grammatical errors have been fixed. The student writes:
\vspace*{-0.3cm}
\begin{quoting}[leftmargin=0.0\parindent]
{\small ``I'm one of your students and I want to tell you something, so I write this e-mail to you.
I am not sure whether you know that most of your students in your class haven't used some programs to solve mathematics problem during middle school and high school. So, most of us don't  know how to use those programs. And you just gave us the assignment to do with those programs. To do that assignment, we have to learn how to use program first. And that's the main point. With no help to learn a program (\textsl{Maxima}, etc.), I think it is too hard and it needs too much time. I have spent more than four or five hours on that assignment and I haven't finished yet. So if you want to give us some assignments, please just make us solve the problems from the book. Now I'm using \textsl{Maxima} to do that assignment because that's the only one I know, and until now I don't know how to use it properly. Please just let us solve the problems from the book if you want to give us homework. Please.\\ P.S. I hope this e-mail will not make you feel bad and all I want is to communicate with you, professor.\\ From: 2016.xxx.xxx (GEDB001-57).'' \par}
\end{quoting}

\vspace*{-0.3cm}
It turns out that this experience was not unique. Other students also encountered difficulty although we asked them to view a video recording explaining a similar problem. After investigation, the problem lies in the software installation, possibly due to bugs in the executable file. Other students were able to complete the assignment correctly with little or no further guidance. We mistakenly presumed that this assignment is manageable, due to the video recording and the fact that our students are `digital natives'~\citep{prensky2001digital}. After this incidence, we did not assign any CAS-related problems during the second half of the semester, but instead focused on problem-solving sessions using pencil and paper. The end of semester feedback obtained also voices a similar concern on {\sl Maxima} operating difficulties. 

\section{Discussion} \label{discussion}

In this study, we have explored students' academic performance and perception after implementing flipped classroom pedagogy in several SVC sections in an EMI--CHC context. In what follows, we discuss some limitations, conclusions and future possibilities and implications. 

\subsection{Limitation}
There are a number of limitations encountered in this study. First, 310 students do not provide a sizable data set. There is also an imbalance between the number of students in the experimental groups (225) and the control group (85). A larger number of students is desirable; if this is achieved, the outcome regarding the effectiveness of flipped classroom pedagogy might be strongly confirmed. Second, there are no organized pre-tests and post-tests in this study, whether for English proficiency, prior mathematical background (basic skill test or Calculus concepts), or CAS experience. 

Third, Type C provides an anomalous result. Since only one topic is flipped, a fair and more precise comparison should be conducted with an identical topic. However, since the common final exam covers this particular topic, we decided to compare the exam score and final grade instead. The student body in Type C is quite diverse as well, with a high number of repeaters and students of international origin. These factors translate to diverse mathematical and language proficiencies. Although the latter does not seem to be posing a significant problem in teaching and learning, the former affects the academic performance significantly.

Fourth, considering that our students are `digital natives', we assumed they might embrace the use of technology in mathematics education. But we discovered that extra guidance on technical aspects of the CAS may be necessary before the students become comfortable using it. Success in this step can translate to other positive outcomes related to the main focus of our pedagogy by encouraging active participation while flipping the class. 

Fifth, one author (NK) had first time recordings and the result is not as professional as those from {\sl Khan Academy} and other online open courseware video lectures. As addressed by~\cite{guo2014video}, in addition to videos of shorter duration, the type of videos that intersperse an instructor's talking-head, with instructors who speak with high enthusiasm and personal feeling are the ones which are more engaging to the students. Although the recordings are imperfect, we expected some feedback from students with limited English proficiency, who might benefit from pausing and re-watching the videos often; we did not receive such feedback.

\subsection{Conclusion}

We observed that a number of factors might hinder the students' learning process. These include a passive-receptive learning style, potential language and cultural barriers, CAS experience, and diverse academic backgrounds. We have implemented and fostered various pedagogical approaches to address these differences. By creating a class atmosphere with active learning and interactive engagement among our students, we hoped to improve not only student-instructor and student-student communication but also overall academic achievement.

The flipped classroom pedagogy implemented in this study is an intersection of three components of theoretical framework: inverted Bloom's taxonomy of educational learning objectives, teaching and learning in an English-medium environment for non-native English speakers of CHC origin and technology adaptation in teaching and learning. We have investigated the effect of implementing the flipped classroom pedagogy on academic performance  and we also gathered students' feedback on the flipped approach they experienced.

A quantitative result from academic performance indicates that there was no significant difference in the final letter grade distribution across the four types of instruction. This may be due relative grading. On the other hand, there was a significant difference in the exam score between Type A and Type C albeit a small practical significance effect size. The exam type may explain this, for example, a collectively prepared and centrally administered common exam or an independently prepared exam. Another possible reason may be that Type C has anomalous demographics with an unusually high number of repeaters and students without a South Korean secondary education, which accounts for more than a third of the class population in each case.

After experiencing the flipped classroom student perceptions were mixed. Generally, positive comments are related to improved communication and appreciation for the instructor's efforts to encourage interactive engagement among students. Negative comments include difficulty understanding material from the video recordings and some hurdles while adapting to technology for Calculus learning, particularly when {\sl Maxima} was introduced and embedded in classroom activities. We may have overestimated the curiosity and adaptability levels of our digitally native students when introducing {\sl Maxima}; which was utilized to provide quality visualizations, improve understanding of Calculus concepts, and enhance interaction and communication. Indeed, to be successful in teaching with any technological tool, we must cultivate curiosity and flexibility in our students so that they may embrace and benefit from its utilization. On the other hand, although we would like to minimize competition and maximize collaboration among the students by having small group activities, some students felt that an element of competition should be maintained. Collaboration and competition are often opposing ideas, thus finding the balance and successfully integrating them while instructing to aid learning requires careful thought and preparation. 

Due to language and cultural barriers, many students were hesitant to express their opinion publicly in front of their peers. Nevertheless, they expressed their concerns and feedback in a written form. This suggests improved student-instructor communication when it comes to teaching and learning process, particularly in the EMI--CHC context. In previous courses, several students expressed difficulty following classes in English. In this study, we did not encounter this particular common concern across the four types of instruction. The participants of the study may possess higher English proficiency than those in previous courses. The concerns emphasized were related to the implementation of the flipped classroom pedagogy, not to EMI related issues. Nevertheless, 40\% to 70\% of participants wrote their feedback in Korean, depending on the instruction types; suggesting that a large percentage of students are more comfortable providing feedback in their native language. On the flip side, the students' use of Korean may not detract from the achievement of the following goals, helping students develop critical thinking skills and creating a class environment with active learning and interactive engagement.

\subsection{Future implication}
There are a number of future implications that can be drawn from this study. Firstly, more rigorous and further pedagogical theory for the flipped classroom is needed to explain why the pedagogy works well for some cases and not as well for the others, such as in the EMI--CHC context. For instance, a `modified' Bloom taxonomy with an additional dimension including language and cultural dimensions may be developed. Could it be that the `traditional' flipped classroom cannot overcome the language (EMI) and cultural (CHC) barriers? Would flipped classroom pedagogy work better if both instructors had higher Korean proficiency levels? Or would an EMI `flipped classroom' be more effective if both instructors and students had high proficiency levels in English when both are non-native English speakers and possess common learning culture?


{\footnotesize
\bibliography{2018EMICalFlipBib} 
\bibliographystyle{apalike}
}

\newpage 
\small \sffamily
\begin{center}
\textbf{Appendix A: Assessment for Flipped Classroom Type A (Spring 2016)}
\end{center}
\begin{multicols}{2}
\begin{center}
\textbf{Midterm Test}\\
\textbf{Tuesday, 26 April 2016 \quad 19:30--20:30}
\end{center}
\begin{enumerate}[leftmargin=1em] \vspace{0.01cm}
\item Use a linear approximation to estimate $\sqrt[3]{1001}$. \hfill {\color{blue} [10]}

\item A water tank has the shape of an inverted circular cone with the base radius of 4~m and height 8~m.
Water is leaking out of the tank at the rate of 4~m$^3/$minute at the same time that water is being pumped into the tank at the rate of 2~m$^3/$minute.
How fast is the water level dropping when the height of the water is 3~m? \hfill {\color{blue} [10]}

\item An observer stands at a point $P$, one unit away from a track.
$A$ and $B$ start at the point $S$ in the figure and run along the track.
$A$ runs four times as fast as $B$. Find the maximum value of the observer's angle of sight $\theta$ between the runners.
\hfill {\color{blue} [15]}

\item Find a function $p(t)$ and a constant $c$ such that \hfill {\color{blue} [10]}
\begin{equation*}
2016 + \int_{c}^{x} \frac{p(t)}{t^2} \, dt = 7x^2.
\end{equation*}

\item Evaluate \hfill {\color{blue} [15]}
\begin{equation*}
\int_{1}^{\infty} \frac{x}{x^3 + 1} \, dx.
\end{equation*}

\item Find the limit \hfill {\color{blue} [15]}
\begin{equation*}
\lim_{x \rightarrow 0} \frac{(1 + x)^{\frac{1}{x}} - e}{x}
\end{equation*}
\textsl{Hint:} 
\begin{equation*}
\lim_{x \rightarrow 0} (1 + x)^{\frac{1}{x}} = e.
\end{equation*}

\item Let $f(x) = x^3 + ax - 1 = 0$. If we have the second and the third approximations $x_2 = \frac{1}{2}$ and $x_3 = b$ to the root of $f(x) = 0$ from Newton's method with initial approximation $x_1 = 1$, find the value of $ab$. \hfill {\color{blue} [10]}

\item Consider a lima\c{c}on with inner loop where the equation is given by the polar curve $r = 1 \, + \, 6 \cos \theta$, $0 \leq \theta \leq 2\pi$.
Find all points in the lima\c{c}on where the tangent lines are horizontal. \hfill {\color{blue} [15]} 
\begin{center}
\begin{tikzpicture}[scale=1.25]
\draw[->] (-1,0) -- (4,0);
\draw[->] (0,-2.2) -- (0,2.3);
\draw node [black] at (3,2) {\scriptsize{$r = 1 + 6 \cos \theta$}};
\draw node [black] at (2.6,-0.3) {\scriptsize{$5$}};
\draw node [black] at (3.6,-0.3) {\scriptsize{$7$}};
\draw node [black] at (4.2,0) {\scriptsize{$x$}};
\draw node [black] at (0,2.5) {\scriptsize{$y$}};
\foreach \x/\xtext in {0.5/1, 1/, 1.5/3, 2/, 2.5/ , 3/ , 3.5/}
    \draw[shift={(\x,0)}] (0pt,2pt) -- (0pt,-2pt) node[below] {\scriptsize{$\xtext$}};
\foreach \y/\ytext in {-2/-4, -1.5/, -1/-2, -0.5/, 0/, 0.5/, 1/2, 1.5/, 2/4}
    \draw[shift={(0,\y)}] (2pt,0pt) -- (-2pt,0pt) node[left] {\scriptsize{$\ytext$}};    
\draw[color=blue,domain=0:6.28,samples=200,smooth] plot (xy polar cs:angle=\x r,radius = {.5+3.0*cos(\x r)});    
\end{tikzpicture}
\end{center}
\end{enumerate}

\columnbreak
\begin{center}
\textbf{Final Exam}\\
\textbf{Tuesday, 21 June 2016 \quad 19:30--20:30}
\end{center}
\begin{enumerate}
\item For the parametric curve $x = e^t - t$, $y = 4 e^{\frac{t}{2}}$, $0 \leq t \leq 3$, let $L(t)$ be the length of the curve from $0$ to $t$.
Find a number $t$ such that ${\displaystyle \frac{dL}{dt} = 4}$. \hfill {\color{blue} [10]}

\item Let ${\displaystyle f(x) = \frac{x^2}{(1 - 2x)^3}}$. \hfill {\color{blue} [10]}
\begin{itemize}
\item[(a)] Find a power series representation for $f(x)$.
\item[(b)] Find $f^{(2016)}(0)$.
\end{itemize}

\item Find the orthogonal trajectories of the family of curves ${\displaystyle y = \frac{k}{x}}$. \hfill {\color{blue} [10]}

\item A rumor tends to spread according to the logistic differential equation ${\displaystyle \frac{dy}{dt} = 0.3 y - 0.0001 y^2}$,
where $y$ is the number of people in the community who have heard the rumor and $t$ is the time in days.
Assume that there were 10 people who knew the rumor at initial time $t = 0$. 
\begin{itemize}[leftmargin=1.75em]
\item[(a)] Find the solution of the differential equation. \hfill {\color{blue} [5]}
\item[(b)] How many days will it take for half of the carrying capacity to hear the rumor? \hfill {\color{blue} [5]}
\item[(c)] What are the equilibrium solutions? \hfill {\color{blue} [5]}
\end{itemize}

\item Evaluate the following infinite sum \hfill {\color{blue} [10]}
\begin{equation*}
\sum_{n = 1}^{\infty} n^2 \frac{2^n - 1}{3^n}.
\end{equation*}

\item Find the radius of convergence of the power series ${\displaystyle \sum_{n = 0}^{\infty} a_x x^n}$, where
$\{ a_n \}$ is defined by $a_{n + 1} = \ln(a_n + 3)$, $n \geq 1$, $a_1 = 1$. \hfill {\color{blue} [15]}

\item Find the radius of convergence and the interval of convergence of  \hfill {\color{blue} [15]}
\begin{equation*}
\sum_{n = 2}^{\infty} \frac{(x - 1)^n}{ n (\ln n)^2}.
\end{equation*}

\item Find the area of the region \textsl{inside} the cardioid $r = 1 - \cos \theta$ and also \textsl{inside} the circle $r = \cos \theta$.
\hfill {\color{blue} [15]}
\end{enumerate}
\end{multicols}

\newpage
\begin{center}
\textbf{Appendix B: Assessment for Flipped Classroom Type B and Traditional Classroom Type~D (Spring~2015)}
\end{center}
\begin{multicols}{2}
\begin{center}
\textbf{Midterm Test}\\ 
\textbf{Tuesday, 21 April 2015 \quad 19:30--20:30}
\end{center}
\begin{enumerate}[leftmargin=1em]
\item  Find an equation for the tangent line to the following curve at $x = 1$: \hfill {\color{blue} [10]}
\begin{equation*}
y = \sqrt[x]{x^{2015} + \ln x}.
\end{equation*}
\item Find a function $f$ and a number $a$ satisfying \hfill {\color{blue} [15]}
\begin{equation*}
1 + \int_{a}^{x} \frac{f(t)}{t^2} \, dt = \sin^{-1} x.
\end{equation*}
\item Consider the equation $(x^2 + y^2)^{3/2} = x^2$ given by the Cartesian coordinate system. \hfill {\color{blue} [15]}
\begin{itemize}[leftmargin=1.75em]
\item[(a)] Convert the given equation to polar equation and sketch the graph, and discuss the symmetry with respect to the origin, $x$-axis, and $y$-axis.
\item[(b)] Let $P$ be the point on the first quadrant where it has the horizontal tangent line.
Find the equation of the line passing two points $P$ and the origin.
\end{itemize}
\item  Consider the set $C$ of points $(x,y)$ satisfying the equation $y^2 = x^3$. 
What is the distance from the point $P(1/2,0)$ to $C$? \hfill {\color{blue} [15]}
\item  Arrange the following three numbers in an increasing order: \hfill {\color{blue} [15]}
\begin{eqnarray*}
A &=& \lim\limits_{x \rightarrow 0^+} x^{\sin x} \\
B &=& \lim\limits_{x \rightarrow 0} \frac{x - \sin x}{\tan x - x} \\
C &=& \lim\limits_{x \rightarrow 0} \frac{e^x - \sin x}{\tan x - x}
\end{eqnarray*}
\item  Let ${\displaystyle F(x) = \int_{0}^{x} \sqrt{t} \sin t \, dt}$. Decide whether $F$ has a local maximum or minimum at the least positive critical number. \hfill {\color{blue} [15]}
\item For what values of $p$ is the integral
\begin{equation*}
\int_{0}^{1} x^p \ln x \, dx
\end{equation*}
convergent? \hfill {\color{blue} [15]}
\end{enumerate}

\columnbreak
\begin{center}
\textbf{Final Exam}\\
\textbf{Tuesday, 16 June 2015 \quad 19:30--20:30}
\end{center}
\begin{enumerate}
\item Find the area of the region in the first quadrant enclosed by the curves $y = 2^{1 - x}$, $y = x^3$, and the $y$-axis. \hfill {\color{blue} [10]}

\item There is a metal ball of radius 1~cm. We want to make a ring by drilling a hole of radius $r$~cm through the center of the ball 
so that the volume of the ring is equal to $\frac{1}{8}$ times the volume of the ball. What is the radius $r$? \hfill {\color{blue} [15]}

\item Determine whether the following series converges. \hfill {\color{blue} [15]}
\begin{itemize}[leftmargin=1.75em]
\item[(a)] ${\displaystyle \sum_{n = 1}^\infty \frac{2^n + 1}{n 2^n + 1}}$
\item[(b)] ${\displaystyle \sum_{n = 1}^\infty \frac{2^n}{1 + (\ln n)^n}}$
\item[(c)] ${\displaystyle \sum_{n = 1}^\infty n^6 e^{-n^{16}}}$
\end{itemize}

\item Find the area of the region that lies inside the curve $r = 2 + \cos 2\theta$ and outside the curve $r = 2 + \sin \theta$. \hfill {\color{blue} [15]}

\item Find the orthogonal trajectories of the family of curves given as ${\displaystyle y = \frac{1 - kx}{kx}}$, where $k$ is a constant. \hfill {\color{blue} [15]}

\item For ${\displaystyle f(x) = \frac{x}{x^2 + 3x + 2}}$, find $f^{(2015)}(0)$. \hfill {\color{blue} [15]}

\item Find the sum of the following series. \hfill {\color{blue} [15]}
\begin{equation*}
\sum_{n = 0}^\infty \frac{(-1)^n}{(2n + 1) 3^n} = 1 - \frac{1}{3 \times 3} + \frac{1}{5 \times 3^2} + \frac{1}{7 \times 3^3} + \dots
\end{equation*}
\end{enumerate}
\end{multicols}

\newpage
\begin{center}
\textbf{Appendix C: Assessment for Flipped Classroom Type C (Fall 2015)}
\end{center}
\begin{multicols}{2}
\begin{center}
\textbf{Midterm Test}\\
\textbf{Wednesday, 21 October 2015 \quad 10:30--11:30}
\end{center}
\begin{enumerate}[leftmargin=1em]
\item Evaluate the following limits, if they exist. Justify your steps.
\begin{itemize}[leftmargin=1.75em]
\item[(a)] ${\displaystyle \lim_{x \rightarrow -2} \frac{x + 2}{x^3 + 8}}$ \hfill {\color{blue} [10]}
\item[(b)] ${\displaystyle \lim_{x \rightarrow -4} \frac{\sqrt{x^2 + 9} - 5}{x + 4}}$ \hfill {\color{blue} [10]}\\
\end{itemize}

\item Consider the function $f(x) = x^4 - 2x^2 + 3$.
\begin{itemize}[leftmargin=1.75em]
\item[(a)] Find the intervals on which $f$ is increasing or decreasing. \hfill {\color{blue} [5]}
\item[(b)] Find the local maximum and minimum values of~$f$. \hfill {\color{blue} [5]}
\item[(c)] Find the intervals of concavity and the inflection points. \hfill {\color{blue} [10]}\\
\end{itemize}

\item 
\begin{itemize}[leftmargin=1.75em]
\item[(a)] Use implicit differentiation to find an equation of the tangent line to the cardioid: $x^2 + y^2 = (2x^2 + 2y^2 - x)^2$  at $\left(0,\frac{1}{2}\right)$. \hfill {\color{blue} [10]}
\item[(b)] Using the logarithmic differentiation method, show that the derivative of the function $y = (\sin x)^{\ln x}$ is given by \hfill {\color{blue} [10]}
\begin{equation*}
\frac{dy}{dx} = (\sin x)^{\ln x} \left(\ln (\sin x)^{1/x} + \ln x^{\cot x}\right).\\
\end{equation*}
\end{itemize}

\item
\begin{itemize}[leftmargin=1.75em]
\item[(a)] The top of a ladder slides down a vertical wall at a rate of 0.15~m/s.
At the moment when the bottom of the ladder is 3~m from the wall, it slides away from the wall at a rate of 0.2~m/s.
How long is the ladder? \hfill {\color{blue} [10]}

\item[(b)] If 1200~cm$^2$ of material is available to make a box with a square base and an open top,
find the largest possible volume of the box. \hfill {\color{blue} [10]}
\end{itemize}

\item
\begin{itemize}[leftmargin=1.75em]
\item[(a)] Using the Fundamental Theorem of Calculus, find the derivative of the function  \hfill {\color{blue} [10]}
\begin{equation*}
y = \int_{0}^{\tan x} \sqrt{t + \sqrt{t}} \, dt.
\end{equation*}

\item[(b)] Using the substitution rule, evaluate the integral \hfill {\color{blue} [10]}
\begin{equation*}
\int \frac{\tan^{-1} x}{1 + x^2} \, dx.\\
\end{equation*}
\end{itemize}
\end{enumerate}

\columnbreak
\begin{center}
\textbf{Final Exam}\\
\textbf{\quad Wednesday, 16 December 2015 \quad 10:30--11:30}
\end{center}
\begin{enumerate} 
\item 
\begin{itemize}[leftmargin=1.75em]
\item[(a)] Using the disk method, find the volume of the solid obtained by rotating the region bounded by $y = \frac{1}{x}$, $x = 1$, $x = 2$ and $y = 0$ about the $x$-axis. \hfill {\color{blue} [10]}
\item[(b)] Using the method of cylindrical shell, find the volume of the solid obtained by rotating the region bounded by $y = \sin (x^2)$, $0 \leq x \leq \sqrt{\pi}$ and $y =0$ about the $y$-axis. \hfill {\color{blue} [10]}
\end{itemize}

\item 
\begin{itemize}[leftmargin=1.75em]
\item[(a)] Find the area of the region that is bounded by the polar curve $r = \sin \theta$, where $\frac{\pi}{3} \leq \theta \leq \frac{2\pi}{3}$. \hfill {\color{blue} [10]}
\item[(b)] Find the exact length of the polar curve $r = 3 \sin \theta$, for $0 \leq \theta \leq \frac{\pi}{3}$. \hfill {\color{blue} [10]}
\end{itemize}

\item 
\begin{itemize}[leftmargin=1.75em]
\item[(a)] Using the method of separation of variables, show that the general solution of the differential equation
\begin{equation*}
\frac{dy}{d\theta} = \frac{e^y \sin^2 \theta}{y \sec \theta}
\end{equation*}
can be written implicitly as $3 (y + 1)e^{-y} + \sin^3\theta = C$, where $C$ is a constant. \hfill {\color{blue} [10]}
\item[(b)] Solve the initial value problem \hfill {\color{blue} [10]}
\begin{equation*}
\frac{du}{dt} = \frac{2t + \sec^2 t}{2 u}, \qquad \qquad u(0) = -5. 
\end{equation*}
\end{itemize}

\item 
\begin{itemize}[leftmargin=1.75em]
\item[(a)] Using the integral test, determine whether the series is convergent or divergent \hfill {\color{blue} [10]}
\begin{equation*}
\sum_{n = 2}^{\infty} \; \frac{1}{n \ln n}.
\end{equation*}
Make sure to verify that the corresponding function satisfies the conditions for the integral test.

\item[(b)] Find the radius of convergence and interval of convergence of the series \hfill {\color{blue} [10]}
\begin{equation*}
\sum_{n = 0}^{\infty} \frac{(x - 2)^n}{n^2 + 1}.
\end{equation*}
\end{itemize}

\item
\begin{itemize}[leftmargin=1.75em]
\item[(a)] Find the Taylor series for $f(x) = e^x$ centered at $x = 7$.  \hfill {\color{blue} [10]}

\item[(b)] Use a binomial series to obtain a Maclaurin series for the function \hfill {\color{blue} [10]}
\begin{equation*}
g(x) = \frac{1}{(1 - x)^2}.
\end{equation*}
\end{itemize}
\end{enumerate}
\end{multicols}
\end{document}